\documentclass[a4paper,11pt]{article}
\usepackage{mfpaperstuff2}
\usepackage{mfabacus}
\usepackage{multirow}
\usetikzlibrary{calc}

\begin{document}

\renewcommand\bottomfraction{0.7}
\newcommand\rem[2]{\operatorname{rem}_{#1}(#2)}
\newcommand\nor[2]{\operatorname{nor}_{#1}(#2)}
\newcommand\blu{\Yfillcolour{blue!25!white}}
\newcommand\wht{\Yfillcolour{white}}
\newcommand\albe{\al\be}
\newcommand\alga{\al\ga}
\newcommand\gaal{\ga\al}
\newcommand\bega{\be\ga}
\newcommand\hand[1]{\operatorname{hand}(#1)}
\newcommand\lel{\operatorname{ll}}
\newcommand{\resi}{\operatorname{res}}
\newcommand{\fx}{x}
\newcommand{\fy}{y}
\newcommand{\rez}{\widehat{\resi}}
\newcommand\sw[1]{#1^\leftrightarrow}
\newcommand\wt[1]{\operatorname{wt}(#1)}
\newcommand\lra\longrightarrow
\newcommand\Hom{\operatorname{Hom}}
\newcommand\tricky{exceptional\xspace}
\newcommand\untricky{unexceptional\xspace}
\newcommand\Tricky{Exceptional\xspace}
\renewcommand\rt[1]{\rotatebox{90}{$#1$}}
\newcommand\bp[2]{\big(#1\ \big|\ #2\big)}
\newcommand\tjs{\stackrel{\mathrm{JS}}\longrightarrow}
\newcommand\js{Jantzen--Schaper\xspace}
\newcommand\ee[1]{\mathrm{e}_{#1}}
\newcommand\ff[1]{\mathrm{f}_{#1}}
\newcommand\eed[2]{\mathrm{e}_{#1}^{(#2)}}
\newcommand\ffd[2]{\mathrm{f}_{#1}^{(r)}}
\newcommand\dm{^\diamond}
\newcommand\bi{\overline i}
\newcommand\bj{\overline j}
\newcommand\bip{bipartition\xspace}
\newcommand\bips{bipartitions\xspace}
\newcommand\zez{\bbz/e\bbz}
\newcommand\xr{nucleus\xspace}
\newcommand\dl[3]{\operatorname{a}_{#2}^{(#3)}{}{\downarrow}^{#1}}
\newcommand\dw[1]{\downarrow^{#1}}
\newcommand\ddu[3]{\downarrow^{#1,#2}\uparrow_{#3}}
\newcommand\spe[1]{\operatorname{S}_{#1}}
\newcommand\wspe[1]{{\operatorname{S}}^\leftrightarrow_{#1}}
\newcommand\hspe[1]{\hat{\operatorname{S}}_{#1}}
\newcommand\hwspe[1]{\hat{\operatorname{S}}^\leftrightarrow_{#1}}
\newcommand\smp[1]{\operatorname{D}_{#1}}
\newcommand\dspe[1]{\operatorname{S}'(#1)}
\newcommand\dsmp[1]{\operatorname{D}'(#1)}
\newcommand\dn[2]{[\spe{#1}\nolinebreak:\nolinebreak\smp{#2}]}
\newcounter{tablecase}
\def\caselabel{xx}
\newcommand\nextcase{\stepcounter{tablecase}\thetablecase.}
\DeclareRobustCommand\rnextcase[1]{\refstepcounter{tablecase}\thetablecase.\label{#1}}

\newcommand\ol[1]{\widebar{#1}}

\newcommand\hhh{\mathcal{H}_}
\newcommand\thhh{\hat{\mathcal{H}}_}

\title{Decomposition numbers for weight $3$ blocks\\of Iwahori--Hecke algebras of type $B$}

\msc{20C08, 20C30, 05E10}

\author{\begin{tabular}{cc}{\Large Matthew Fayers}&{\Large Lorenzo Putignano}\\[3pt]
{\normalsize Queen Mary University of London}&{\normalsize Universit\`a degli Studi di Firenze}\\
{\normalsize\texttt{\normalsize m.fayers@qmul.ac.uk}}&{\normalsize\texttt{\normalsize lorenzo.putignano@unifi.it}}
\end{tabular}}

\renewcommand\auth{Matthew Fayers \& Lorenzo Putignano}
\runninghead{Weight $3$ blocks of Iwahori--Hecke algebras of type $B$}

\toptitle

\begin{abstract}
  Let $B$ be a weight-$3$ block of an Iwahori--Hecke algebra of type $B$ over any field. We develop the combinatorics of $B$ to prove that the decomposition numbers for $B$ are all $0$ or $1$.
\end{abstract}

\tableofcontents

\section{Introduction}

In the representation theory of the symmetric groups, a fundamental problem is to determine the \emph{decomposition numbers} in characteristic $p$: the composition multiplicities of $p$-modular simple representations in $p$-modular reductions of ordinary simple representations. A complete solution to this problem seems far out of reach, but solutions are known in various special cases. One of these results (due to the first author \cite{mfwt3}) says that the decomposition numbers in a $p$-block of weight $3$ in characteristic $p\gs5$ are all either $0$ or $1$; this in particular allows these decomposition numbers to be computed quickly using the Jantzen--Schaper formula.

The decomposition number problem naturally extends to the Iwahori--Hecke algebras of type $A$, where we seek the composition multiplicities of the simple modules in the Specht modules. Many of the known results for symmetric groups have been proved for the Iwahori--Hecke algebras as well, including the result for blocks of weight $3$. More recently, this problem has been extended to the Iwahori--Hecke algebras $\hhh n$ of type $B$. These were first systematically studied by Dipper, James and Murphy \cite{djm}, who introduced appropriate definitions of Specht modules labelled by \bips of $n$. This naturally leads to a combinatorial approach to the corresponding decomposition number problem, underpinned by the combinatorics of \bips. The first author's combinatorial definition of the weight of a \bip, and hence the weight of a block \cite{mfweight}, allows us to approach this problem systematically by addressing blocks of a fixed weight. As in type $A$, the blocks of weight $0$ are precisely the simple blocks, and blocks of weight $1$ admit a very simple description (analogous to symmetric group blocks with cyclic defect groups). The first author \cite{mftypebwt2} gave a combinatorial formula for decomposition number for blocks of weight $2$, analogous to Richards's formula \cite{rich} for type~$A$.

In this paper we address blocks of weight $3$. Our main result is as follows.

\begin{thm}\label{main}
Suppose $\bbf$ is a field, $q\in\bbf$ is non-zero, and $k_1,k_2\in\bbz$. Let $\hhh n$ denote the Iwahori--Hecke algebra of type $B$ over $\bbf$ with parameters $q,q^{k_1},q^{k_2}$, and let $B$ be a block of $\hhh n$ of weight $3$. Then the decomposition numbers for $B$ are all either $0$ or $1$.
\end{thm}

Note in particular that (unlike in type $A$) \cref{main} applies even when $\bbf$ has characteristic $2$ or $3$. We conjecture that this result extends to all weight $3$ blocks of \emph{Ariki--Koike algebras}; these are higher-level Hecke algebras generalising the Iwahori--Hecke algebras of types $A$ and $B$. A particular case of this conjecture (the case of `tree blocks') has already been proved by Lyle and Ruff in \cite{lyru}.

The proof of \cref{main} is by induction on $n$, using the Brundan--Kleshchev branching rules which govern induction and restriction between $\hhh n$ and $\hhh m$, for $m<n$. An important part of our work is a combinatorial description of all the \bips labelling Specht modules in a given block, and an analysis of which Specht modules $\spe\la$ cannot be shown to be multiplicity-free via the inductive hypothesis. For the most difficult cases, we need to use the cyclotomic \js formula due to James and Mathas \cite{jmjs}.

\begin{acks}
Most this work was developed while the second author was visiting Queen Mary University of London between October 2022 and April 2023. The second author expresses his sincere gratitude to the University of Florence for having funded the visit and to the first author and the whole QMUL maths department for their hearty hospitality.
\end{acks}

\section{Background}

In this section we summarise some of the background results we shall need. For notation we mostly follow the first author's paper \cite{mftypebwt2}, where more details and motivation may be found.

\subsection{Basic notation}

In this paper $\bbn$ denotes the set of positive integers. Throughout, we fix an integer $e\gs2$. For any integer $a$, we write $\ol a=a+e\bbz=\lset{a+em}{m\in\bbz}$, and $\zez=\lset{\ol a}{a\in\bbz}$. The set $\zez$ admits an additive action of $\bbz$ via $a+\ol b=\widebar{a+b}$.

In general, given any tuple $\lset{x_i}{i\in\zez}$ of objects indexed by $\zez$, and given an integer $a$, we may write $x_a$ to mean $x_{\ol a}$. Similarly, given any terminology involving an element of $\zez$ (such as `$i$-node'), we may substitute an integer $a$ in place of $\ol a$ (so we may define an $a$-node to mean an $\ol a$-node).

Throughout the paper, we fix a pair $\ka=(\ka_1,\ka_2)\in(\zez)^2$. A pair of integers $k_1,k_2$ such that $\ka_1=\ol{k_1}$ and $\ka_2=\ol{k_2}$ is called a \emph{bicharge} for $\ka$.

If $M$ and $N$ are modules for an algebra, then we write $M\sim N$ to mean that $M$ and $N$ have the same composition factors with multiplicity.

\subsection{Partitions and \bips}

A \emph{partition} is an infinite weakly-decreasing sequence $\la=(\la_1,\la_2,\dots)$ of integers which are eventually zero. We write $\card\la$ for the sum $\la_1+\la_2+\dots$, and say that $\la$ is a partition of $\card\la$. When writing partitions, we omit trailing zeroes and group together equal parts with a superscript. We write $\vn$ for the unique partition of $0$. For any partition $\la$, we write $h(\la)$ for the number of non-zero parts of $\la$, which we call the \emph{length} of $\la$.

A \emph{\bip} is an ordered pair $\la=\bp{\la^{(1)}}{\la^{(2)}}$ of partitions, which we call the \emph{components} of $\la$. We write $\card\la=\card{\la^{(1)}}+\card{\la^{(2)}}$, and say that $\la$ is a \bip of $\card\la$.

The Young diagram of a \bip $\la$ is the set
\[
[\la]=\lset{(r,c,a)\in\bbn^2\times\{1,2\}}{c\ls\la^{(a)}_r},
\]
whose elements we call the \emph{nodes} of $\la$. In general, a node means an element of $\bbn^2\times\{1,2\}$. A node $(r,c,a)$ of $\la$ is \emph{removable} if it can be removed from $[\la]$  to leave the Young diagram of a smaller \bip (that is, if $c=\la^{(a)}_r>\la^{(a)}_{r+1}$), while a node not in $\la$ is an addable node of $\la$ if it can be added to $[\la]$ to leave a Young diagram.

We depict $[\la]$ by drawing the Young diagram of $\la^{(1)}$ above the Young diagram of $\la^{(2)}$ (using the English convention). Accordingly, we will say that a node $(r,c,a)$ is \emph{above} or \emph{higher than} a node $(s,d,b)$ if either $a<b$ or $a=b$ and $r<s$.

If $\la$ and $\mu$ are two \bips of $n$, then we say that $\la$ \emph{dominates} $\mu$ (and write $\la\dom\mu$) if
\[
\la^{(1)}_1+\dots+\la^{(1)}_r\gs\mu^{(1)}_1+\dots+\mu^{(1)}_r\qquad\text{and}\qquad\card{\la^{(1)}}+\la^{(2)}_1+\dots+\la^{(2)}_r\gs\card{\mu^{(1)}}+\mu^{(2)}_1+\dots+\mu^{(2)}_r
\]
for every $r\gs1$. Another way of saying this is that $[\mu]$ can be obtained from $[\la]$ by moving one or more nodes to lower positions.

If $\la$ is a partition, the \emph{conjugate} partition $\la'$ is defined by
\[
(\la')_r=\card{\lset{c\in\bbn}{\la_c\gs r}}.
\]
If $\la$ is a \bip, the conjugate \bip is $\la'=\bp{{\la^{(2)}}'}{{\la^{(1)}}'}$. In other words, $\la'$ is obtained by interchanging the two components, and then replacing each component with the conjugate partition. Note then that if $\la,\mu$ are two \bips, then $\la\dom\mu$ \iff $\mu'\dom\la'$.

Now suppose we have fixed a pair $\ka=(\ka_1,\ka_2)\in(\zez)^2$. The \emph{residue} of a node $(r,c,a)$, denoted $\resi(r,c,a)$, is $c-r+\ka_a$. If $i\in\zez$, then an \emph{$i$-node} means a node of residue $i$. If $\la$ is a \bip, then the \emph{content} of $\la$ is the multiset of the residues of the nodes of $\la$. For any $i\in\zez$, we write $\rem i\la$ for the number of removable $i$-nodes of $\la$.

We will also need to consider rim hooks. The \emph{rim} of a \bip $\la$ is the set of nodes $(r,c,a)\in[\la]$ such that $(r+1,c+1,a)\notin[\la]$. A \emph{rim hook} of $\la$ is a connected subset of the rim which can be removed to leave a smaller \bip. If $H$ is a rim hook, the \emph{hand node} $\hand H$ is the top-rightmost node of $H$, and the \emph{leg length} $\lel(H)$ is the difference in height between the highest and lowest nodes in $H$.

\subsection{Beta-sets and the abacus}

We often use James's abacus for partition combinatorics. With $e$ fixed as above, we take an abacus with $e$ vertical runners labelled $0,\dots,e-1$ from left to right, and with positions marked $0,1,\dots$ from left to right along successive rows of the abacus, working down the page. Now given a partition $\la$ and an integer $k\gs h(\la)$, we define the \emph{beta-set for $\la$} with charge $k$ to be the set $\lset{\la_r+k-r}{1\ls r\ls k}$. The \emph{\abd} for $\la$ with charge $k$ is constructed by placing a bead on the abacus at position $\la_r+k-r$ for each $1\ls r\ls k$.

Given a \bip $\la=\bp{\la^{(1)}}{\la^{(2)}}$ and a bicharge $(k_1,k_2)\in\bbn^2$ with $k_a\gs h(\la^{(a)})$ for $a=1,2$, the \abd for $\la$ with bicharge $(k_1,k_2)$ is obtained by placing the \abd for $\la^{(1)}$ with charge $k_1$ above the \abd for $\la^{(2)}$ with charge $k_2$.

The abacus is good for visualising rim hooks (and in particular, removable nodes). Given an \abd for a \bip $\la$, removable nodes in component $a$ correspond to positions $x>0$ on the abacus such that in component $a$ there is a bead at position $x$ but no bead at position $x-1$. Given such a $x$, the residue of the corresponding removable node is $x+e\bbz$. Similarly, addable nodes correspond to positions $x$ such that there is a bead at position $x-1$ but not at position $x$.

More generally, rim hooks in component $a$ of $\la$ correspond to pairs of positions $x,y$ where $x>y$ and there is a bead at position $x$ but not at position $y$ in component $a$. The residue of the hand node of the corresponding rim hook is $x+e\bbz$, the leg length is the number of beads between positions $x$ and $y$ (not including position $x$), and removing the rim hook corresponds to moving the bead at position $x$ to position $y$.

\subsection{Restricted \bips}\label{kleshbipsec}

In \cite{arma}, Ariki and Mathas introduced a set of multipartitions to provide a combinatorial labelling for simple modules of Ariki--Koike algebras. They called these \emph{Kleshchev multipartitions}, but here we prefer the term \emph{restricted multipartitions} introduced by Brundan--Kleshchev \cite{bk}. We now recall the definition here, specialising to the case of \bips. We will use the original recursive definition of Ariki--Mathas, though we note that a more efficient definition was found by Ariki, Kreiman and Tsuchioka \cite{akt}.

The definition of restricted \bips depends on our fixed $\ka\in(\zez)^2$. Take a \bip $\la$ and a residue $i\in\zez$. Define the \emph{$i$-signature} of $\la$ to be the sequence of signs obtained by reading the Young diagram of $\la$ from top to bottom, writing $+$ for each addable $i$-node and $-$ for each removable $i$-node. The \emph{reduced $i$-signature} is the subsequence obtained by repeatedly deleting adjacent pairs $-+$. The removable nodes corresponding to the $-$ signs are called the \emph{normal} $i$-nodes of $\la$, and the highest of these (if there are any) is the \emph{good} $i$-node. The addable nodes corresponding to the $+$ signs in the reduced $i$-signature are the \emph{conormal} $i$-nodes of $\la$, and the lowest of these (if there are any) is the \emph{cogood} $i$-node.

For any \bip $\mu$ and any $i\in\zez$, we write $\nor i\mu$ for the number of normal $i$-nodes of $\mu$.

\begin{eg}
Suppose $e=3$, $\ka=(\ol0,\ol1)$ and $\la=\bp{3,2,1^2}{2^3}$. The Young diagram of $\la$ is shown below, with the residues of the addable and removable nodes indicated.
\[
\gyoung(;;;2:0,;;0:1,;:2,0,:2,,;;:0,;;,;;0,:1)
\]
If we let $i=\ol0$, then we see that the $i$-signature of $\la$ is $+--+-$. So the reduced $i$-signature is $+--$. So there are two normal $i$-nodes $(2,2,1)$ and $(3,2,2)$, with $(2,2,1)$ being the good $i$-node. There is a unique conormal (and therefore cogood) $i$-node $(1,4,1)$.
\end{eg}

Now given two \bips $\la$ and $\mu$, write $\mu\stackrel i\longrightarrow\la$ if $\mu$ is obtained from $\la$ by removing the good $i$-node for $i\in\zez$, or equivalently if $\la$ is obtained from $\mu$ by adding the cogood $i$-node. We write $\mu\longrightarrow\la$ to mean that $\mu\stackrel i\longrightarrow\la$ for some $i$, and let $\longleftrightarrow$ be the equivalence relation generated by $\longrightarrow$. Now define the set of restricted \bips (with respect to $\ka$) to be the \bips in the $\longleftrightarrow$-class containing $\bp\vn\vn$.

The definition of restricted \bips derives from the theory of crystals for highest-weight representations of quantum groups. It follows from this theory (in particular, the uniqueness of highest-weight vectors) that the only restricted \bip with no normal nodes of any residue is $\bp\vn\vn$. So to test whether a \bip $\la$ is restricted, one can repeatedly remove arbitrarily-chosen good nodes until there are no good nodes remaining; then $\la$ is restricted \iff the resulting \bip is $\bp\vn\vn$.

We also make an observation which will be helpful later on. Say that a partition $\la$ is \emph{$e$-restricted} if $\la_r-\la_{r+1}<e$ for all $r$. Now the following is a well-known result; it reflects that restricted \bips label the vertices of a level-$2$ highest-weight crystal obtained as a component of the tensor product of two level-$1$ crystals, which themselves have vertices labelled by $e$-restricted partitions.

\begin{lemma}\label{resterest}
Suppose $\la$ is a restricted \bip. Then both $\la^{(1)}$ and $\la^{(2)}$ are $e$-restricted partitions.
\end{lemma}

\begin{pf}
It is an easy exercise to check that if $\la,\mu$ are \bips with $\mu\longrightarrow\la$, and $a\in\{1,2\}$, then $\la^{(a)}$ is $e$-restricted \iff $\mu^{(a)}$ is. Since the components of the empty \bip $\bp\varnothing\varnothing$ are both $e$-restricted, the result follows.
\end{pf}

We will also need to consider \emph{regular} \bips. These are described in \cite[1.3.3]{mftypebwt2}, where they are called \emph{conjugate Kleshchev} \bips. We say that a \bip $\la$ is \emph{regular} if $\la'$ is restricted. In fact we can modify the above definition of restricted \bips to give an alternative description of regular \bips, which will be useful later. To do this, define the $i$-signature of a \bip as above, and define the \emph{antireduced} $i$-signature by successively deleting pairs $+-$. The removable nodes corresponding to $-$ signs in the antireduced $i$-signature are called the \emph{antinormal} $i$-nodes of $\la$, and the lowest of these (if there are any) is the \emph{antigood} $i$-node. Similarly, the addable nodes corresponding to $+$ signs are the \emph{anticonormal} $i$-nodes of $\la$, and the highest of these (if there are any) is the \emph{anticogood} $i$-node. Given two \bips $\la$ and $\mu$, we write $\mu\stackrel i\Longrightarrow\la$ if $\mu$ is obtained from $\la$ by removing the antigood $i$-node for some $i$. Equivalently, if we write $i'=\ka_1+\ka_2-i$, then $\mu\stackrel i\Longrightarrow\la$ \iff $\mu'\stackrel{i'}\longrightarrow\la'$. Now let $\Longrightarrow$ be the union of the relations $\stackrel i\Longrightarrow$, and let $\Longleftrightarrow$ be the equivalence relation generated by $\Longrightarrow$. The \bips in the $\Longleftrightarrow$-class containing $\bp\vn\vn$ are the regular \bips.

\subsection{Iwahori--Hecke algebras}

We fix a field $\bbf$ with $q,Q_1,Q_2$ non-zero elements of $\bbf$, and let $\hhh n$ denote the Iwahori--Hecke algebra of type $B$ with parameters $q,Q_1,Q_2$. This is the unital associative $\bbf$-algebra with generators $T_0,\dots,T_{n-1}$, subject to relations
\begin{alignat*}2
(T_0-Q_1)(T_0-Q_2)&=0,&&
\\
(T_i-q)(T_i+1)&=0&\qquad&\text{for }1\ls i<n,
\\
T_iT_j&=T_jT_i&&\text{if }j-i>1,
\\
T_0T_1T_0T_1&=T_1T_0T_1T_0,
\\
T_iT_{i+1}T_i&=T_{i+1}T_iT_{i+1}&&\text{for }1\ls i<n-1.
\end{alignat*}
The definition of $\hhh n$ is unaffected if we interchange $Q_1$ and $Q_2$, but it will be important for us to consider $(Q_1,Q_2)$ as an ordered pair; in particular, the construction of Specht modules depends on this ordered pair. (In \cref{conjswsec} below, we consider in more detail the effect of interchanging $Q_1$ and~$Q_2$.)

As explained in \cite[1.3.2]{mftypebwt2}, using a Morita equivalence result of Dipper--James \cite{dj}, we can restrict attention to the case where $Q_1$ and $Q_2$ are powers of $q$. We also assume that $q$ is a primitive $e$th root of unity, where $e\gs2$; the case where $q$ is not a root of unity works in the same way as the case where $e>n$ (and in fact most of the particular cases of blocks that we consider only arise if $q$ is an $e$th root of unity for $e\ls n$).

So with our fixed pair $\ka=(\ka_1,\ka_2)\in\zez^2$, we will assume henceforth that $Q_1=q^{k_1}$ and $Q_2=q^{k_2}$, where $(k_1,k_2)$ is any bicharge for $\ka$.

The representation theory of $\hhh n$ is based on the theory of \emph{Specht modules}. For each \bip $\la$ of $n$, we write $\spe\la$ for the corresponding Specht module, as defined in Mathas's survey article \cite{mathsurv}. The Specht modules are precisely the cell modules with respect to a particular cellular basis of $\hhh n$, and a central problem in the representation theory of $\hhh n$ is to determine their composition factors. To do this, we need a classification of the simple $\hhh n$-modules. If $\la$ is a restricted \bip, then $\spe\la$ has a unique simple quotient $\smp\la$, and the simple modules $\smp\la$ arising in this way give a complete irredundant set of simple $\hhh n$-modules. Thus the problem of determining the composition factors of the Specht modules amounts to determining the \emph{decomposition numbers} $\dn\la\mu$, where $\la,\mu$ are \bips of $n$ with $\mu$ restricted.

A fundamental result on decomposition numbers is the following.

\begin{propn}[\xcite{mftypebwt2}{Theorem 1.2}]\label{basicdecomp}
If $\la$ and $\mu$ are \bips of $n$ with $\mu$ restricted, then $\dn\mu\mu=1$, while if $\dn\la\mu>0$ then $\la\dom\mu$.
\end{propn}

\subsection{Blocks and weight}

An important tool for studying the decomposition number problem is the classification of blocks, since $\dn\la\mu=0$ unless $\spe\la$ and $\smp\mu$ lie in the same block. It is a consequence of the general theory of cellular algebras that a Specht module always belongs to a single block, and we will abuse notation by saying that $\la$ lies in a block $B$ to mean that $\spe\la$ lies in $B$. Of course, if $\la$ is restricted this means that $\smp\la$ lies in this block too, so to classify the blocks of $\hhh n$ it is sufficient to determine when two Specht modules lie in the same block. The block classification for $\hhh n$ is a special case of the result for Ariki--Koike algebras proved by Lyle and Mathas.

\begin{propn}[\xcite{lyma}{Theorem 2.11}]\label{blockclass}
Suppose $\la$ and $\mu$ are \bips of $n$. Then $\spe\la$ and $\spe\mu$ lie in the same block of $\hhh n$ \iff $\la$ and $\mu$ have the same content.
\end{propn}

As explained in \cite[3.2]{mfweight}, this result can also be formulated in terms of the abacus: if we choose a bicharge $(k_1,k_2)$ for $\ka$ with $k_1,k_2$ sufficiently large, and construct the \abds for $\la$ and $\mu$ with this bicharge, then $\la$ and $\mu$ have the same content \iff for each $i$ the number of beads on runner $i$ (across both components) is the same in each \abd. This leads to another useful way to determine whether two \bips belong to the same block. If $\la$ is a \bip, then for each $i\in\zez$ we define $\de_i(\la)$ to be the number of removable $i$-nodes of $\la$ minus the number of addable $i$-nodes of $\la$. Now we have the following result.

\begin{propn}[\xcite{mfweight}{Proposition 3.2}]\label{samede}
Suppose $\la$ and $\mu$ are \bips of $n$. Then $\la$ and $\mu$ lie in the same block \iff $\de_i(\la)=\de_i(\mu)$ for all $i\in\zez$.
\end{propn}

As a consequence we can define $\de_i(B)$ for a block $B$, meaning $\de_i(\la)$ for any \bip $\la$ in $B$. If $a\in\bbz$, then we may write $\de_a(B)$ to mean $\de_{a+e\bbz}(B)$.

In \cite{mfweight}, the first author introduced the notion of the \emph{weight} $\wt\la$ of a \bip $\la$. Since the weight of $\la$ depends only on the content of $\la$, it is the same for all $\la$ in a given block, and therefore one may define the weight of a block $B$ to be the weight of any \bip in $B$. The weight of $B$ is a non-negative integer which measures the complexity of $B$, and a fruitful approach to the representation theory of $\hhh n$ is to consider blocks of a given small weight.

We will not give the definition of weight in terms of content here, but we will use the results from \cite[1.3.5]{mftypebwt2} showing how to calculate the weight of a \bip using the abacus. We will introduce the necessary notation in \cref{class3sec}.

We will also use the following lemma several times.

\begin{lemma}[\xcite{mfweight}{Lemma 3.6}]\label{36lem}
Suppose $\la$ is a \bip, and that $\mu$ is a \bip obtained by removing $u$ removable $i$-nodes from $\la$. Then
\[
\wt\mu=\wt\la+u(\de_i(\la)-u).
\]
\end{lemma}

This paper is motivated by the following result.

\begin{thm}\label{w012}
Suppose $B$ is a block of $\hhh n$ of weight $0$, $1$ or $2$, and $\la,\mu$ are \bips in $B$ with $\mu$ restricted. Then $\dn\la\mu\ls1$.
\end{thm}

\begin{pf}
The result for blocks of weight $0$ follows from \cite[Theorem 4.1]{mfweight}, and the result for weight $1$ is \cite[Theorem 4.2]{mfweight}. For blocks of weight $2$, the result is given in \cite[Theorems 3.18 \& 4.13]{mftypebwt2}.
\end{pf}

Our aim in this paper is to prove the same result for blocks of weight $3$. In fact for blocks of weight at most $2$ explicit formul\ae{} for the decomposition numbers are known; for weight $3$ we just show that the decomposition numbers are at most $1$ (which is enough to enable the decomposition numbers to be computed quickly using the \js formula). This is analogous to the situation for Hecke algebras of type $A$ \cite{mfwt3}.

\subsection{Conjugation and component-switching}\label{conjswsec}

The set of \bips of $n$ has very natural symmetries, given by conjugation or by simply interchanging components. These symmetries provide useful tools for studying decomposition numbers, which we describe here.

We begin with conjugation, which relates to automorphisms. Let $\theta$ be the automorphism defined  defined by mapping $T_0\mapsto q^{k_1+k_2}T_0^{-1}$ and $T_i\mapsto-qT_i^{-1}$ for $i\gs1$. For any $\hhh n$-module $M$, define the module $M^\theta$ by twisting the $\hhh n$-action by $\theta$. Obviously if $M$ is simple then so is $M^\theta$, and so we can define a bijection $\mu\mapsto\mu\dm$ from the set of restricted \bips to the set of regular \bips by $(\smp\mu)^\theta\cong\smp{(\mu\dm)'}$. Understanding what the automorphism $\theta$ does to Specht modules gives us more information on decomposition numbers, as follows.

\begin{propn}\label{conjdual}
Suppose $\la$ and $\mu$ are \bips of $n$ with $\mu$ restricted. Then $\dn\la\mu=\dn{\la'}{(\mu\dm)'}$.
\end{propn}

\begin{pf}
Our argument is based on \cite[Sections 4 and 5]{mathastilt}, where Mathas recalls the \emph{dual Specht module} $\dspe\la$ for each $\la$. The dual Specht modules are related to the usual Specht modules in two ways. On the one hand, one can check that $\dspe\la$ is isomorphic to $(\spe\la)^\theta$. To see this, we recall how $\spe\la$ is constructed: a certain element $m_\la\in\hhh n$ is defined, and $\spe\la$ is defined to be the right ideal generated by $m_\la$, modulo its intersection with the two-sided ideal generated by the elements $m_\nu$ for $\nu\doms\la$. The dual Specht module $\dspe\la$ is defined similarly using a different element $n_\la$. But one can check that $\theta(m_\la)$ equals $n_\la$ times an invertible element of $\hhh n$, so that the (right or two-sided) ideal generated by $\theta(m_\la)$ coincides with the (right or two-sided) ideal generated by $n_\la$. As a consequence, we obtain $\dn\la\mu=[\dspe\la:(\smp\mu)^\theta]$ for all $\la,\mu$.

On the other hand, Mathas considers contragredient duality $M\mapsto M^\circledast$ induced by the antiautomorphism of $\hhh n$ that fixes each $T_i$. In \cite[Corollary 5.7]{mathastilt} Mathas shows that $\dspe\la\cong(\spe{\la'})^\circledast$. Since contragredient duality preserves composition multiplicities, and the simple modules are self-dual (which follows from the cellularity of $\hhh n$), we find that $[\dspe\la:\smp\mu]=\dn{\la'}\mu$ for every $\la,\mu$. The result follows.
\end{pf}

\begin{cory}\label{basicdecomp2}
Suppose $\la$ and $\mu$ are \bips of $n$, with $\mu$ restricted. Then $\dn{\mu\dm}\mu=1$, while if $\dn\la\mu>0$ then $\mu\dm\dom\la$.
\end{cory}

\begin{pf}
This follows from \cref{basicdecomp,conjdual} using the fact that conjugation reverses the dominance order on \bips.
\end{pf}

Later in the paper we will need to be able to calculate the bijection ${}\dm$ combinatorially. Recall the relations $\stackrel i\longrightarrow$ and $\stackrel i\Longrightarrow$ from \cref{kleshbipsec}.

\begin{propn}\label{dmbij}
Suppose $\la$ and $\mu$ are restricted \bips. Then $\mu\stackrel i\longrightarrow\la$ \iff $\mu\dm\stackrel i\Longrightarrow\la\dm$.
\end{propn}

\begin{pf}

The Brundan--Kleshchev modular branching rules show that the relation $\stackrel i\longrightarrow$ can be described algebraically via
\[
\mu\stackrel i\lra\la\qquad\text{\iff}\qquad\Hom_{\hhh{n-1}}(\smp\mu,\ee i\smp\la)\neq0,             
\]
where $\ee i$ is the $i$th restriction functor (see \cref{branchsec} below). The definition of $\ee i$ and the automorphism $\theta$ shows that $(\ee iM)^\theta=\ee{i'}(M^\theta)$ for any module $M$, 
where as above $i'=\ka_1+\ka_2-i$. Hence
\[
\Hom_{\hhh{n-1}}(\smp\mu,\ee i\smp\la)\neq0\qquad\text{\iff}\qquad\Hom_{\hhh{n-1}}(\smp{(\mu\dm)'},\ee{i'}\smp{(\la\dm)'})\neq0,
\]
and therefore
\[
\mu\stackrel i\lra\la\qquad\text{\iff}\qquad(\mu\dm)'\stackrel{i'}\lra(\la\dm)',
\]
and from \cref{kleshbipsec} this is equivalent to $\mu\dm\stackrel i\Longrightarrow\la\dm$.
\end{pf}

Now we come to component-switching, where we don't need to say as much. If $\la$ is a \bip, let $\sw\la$ denote the \bip $\bp{\la^{(2)}}{\la^{(1)}}$ obtained by switching the two components of $\la$. Note that although the isomorphism type of $\hhh n$ only requires the unordered pair $(\ka_1,\ka_2)$, the definition of the Specht module (and the associated combinatorics of residues, blocks, restricted \bips etc) depends on the ordered pair $(\ka_1,\ka_2)$. We let $\wspe\la$ denote the Specht module constructed using the \bip $\sw\la$ and the ordered pair $(\ka_2,\ka_1)$.

\begin{lemma}\label{wspelen}
Suppose $\la$ is a \bip of $n$. Then $\spe\la$ and $\wspe\la$ have the same composition factors (with multiplicity).
\end{lemma}

\begin{pf}
This uses a specialisation argument. Let $\hat q,\hat Q_1,\hat Q_2$ be algebraically independent indeterminates over $\bbf$, and consider the Iwahori--Hecke algebra $\thhh n$ over $\bbf(\hat q,\hat Q_1,\hat Q_2)$ with parameters $\hat q,\hat Q_1,\hat Q_2$. Then $\thhh n$ is semisimple by \cite[Main Theorem]{arikiss}, with the Specht modules being the simple modules. If we define the Specht module $\hspe\la$ for $\thhh n$ using the \bip $\la$ and the parameters $\hat q,\hat Q_1,\hat Q_2$, and we define the Specht module $\hwspe\la$ using the \bip $\sw\la$ and the parameters $\hat q,\hat Q_2,\hat Q_1$, then a consideration of characters shows that $\hspe\la$ and $\hwspe\la$ are isomorphic. The $\hhh n$-modules $\spe\la$ and $\wspe\la$ are obtained by specialising $\hat q$ to $q$, $\hat Q_1$ to $q^{k_1}$ and $\hat Q_2$ to $q^{k_2}$, where $\ka_1=\ol{k_1}$, $\ka_2=\ol{k_2}$. It is a standard result on specialisation that specialisations of isomorphic modules have the same composition factors.
\end{pf}

\cref{wspelen} allows us to reduce the work we do in considering different types of blocks: given a block $B$ of $\hhh n$, we can study the same block by replacing $\ka$ with $(\ka_2,\ka_1)$ and replacing each \bip $\la$ in $B$ with $\sw\la$. The truth of \cref{main} for $B$ is then independent of which way we view $B$.

\subsection{Branching rules}\label{branchsec}

Now we summarise some essential background on the Brundan--Kleshchev branching rules; the most comprehensive reference for these is \cite{bk}. Recall that the content of a \bip is the multiset of the residues of its nodes, and that two \bips lie in the same block \iff they have the same content. So we may define the content $\calc_B$ of a block $B$ to be the content of any \bip in $B$.

Now suppose $M$ is an $\hhh n$-module lying in a block $B$, and $i\in\zez$. If there is a block $A$ of $\hhh{n-1}$ such that $\calc_A$ is obtained by removing a copy of $i$ from $\calc_B$, then we define $\ee iM$ to be the block component of $M\downarrow_{\hhh{n-1}}$ lying in $A$; otherwise, we set $\ee iM=0$.

The functors $\ee i$ are defined for all $n$, so we can define powers $\ee i^r$. In fact, it is possible to define \emph{divided powers}: for any $\hhh n$-module $M$ and for any $r$ there is an $\hhh{n-r}$-module $\eed irM$ such that $\ee i^rM\cong(\eed irM)^{\oplus r!}$. For any non-zero module $M$, we set $\ep_iM=\max\lset{r}{\ee i^rM\neq0}$.

The following results are essential tools in the study of decomposition numbers. Recall that we write $\nor i\mu$ for the number of normal $i$-nodes of a \bip $\mu$ and $\rem i\la$ for the number of removable $i$-nodes of $\la$.

\begin{propn}\label{brchrules}\indent
\begin{enumerate}
\vst
\item
Suppose $\la$ is a \bip of $n$. Then $\eed ir\spe\la$ has a filtration in which the factors are the Specht modules $\spe\nu$, for all \bips $\nu$ that can be obtained from $\la$ by removing $r$ removable $i$-nodes. In particular, $\ep_i\spe\la=\rem i\la$, and $\eed i{\ep_i\spe\la}\spe\la$ is isomorphic to the Specht module $\spe{\la^-}$, where $\la$ is the \bip obtained by removing all the removable $i$-nodes from $\la$.
\item
Suppose $\mu$ is a restricted \bip of $n$. Then $\ep_i\smp\mu=\nor i\mu$ and $\eed i{\ep_i\smp\mu}\smp\mu\cong\smp{\mu^-}$, where $\mu^-$ is the \bip obtained by removing all the normal $i$-nodes from $\mu$.
\end{enumerate}
\end{propn}

\subsection{The cyclotomic \js formula}\label{jssec}

The cyclotomic \js formula introduced by James and Mathas \cite{jmjs} will be a very useful tool. We describe a special case here, specialising to the case of Hecke algebras of type $B$.

First we note that because every field is a splitting field for $\hhh n$, we can extend the field $\bbf$ without affecting the representation theory. So we assume that the subfield of $\bbf$ generated by $q$ is a proper subfield, and we fix an element $x$ of $\bbf$ lying outside this subfield.

Let $\hat q$ be an indeterminate over $\bbf$, and let $R=\bbf[\hat q^{\pm1}]$. Choose a bicharge $(k_1,k_2)$ for $\ka$ for which $|k_1-k_2|>n$. Let $\fkp$ be the prime ideal in $R$ generated by $\hat q-q$. We want to consider the Iwahori--Hecke algebra of type $B$ over $R$ with parameters $\hat q,\hat Q_1,\hat Q_2$, where
\[
\hat Q_1=\hat q^{k_1},\qquad\hat Q_2=\hat q^{k_2}+x(\hat q-q).
\]

Now given a node $(r,c,a)\in\bbn^2\times\{1,2\}$ define
\[
\rez(r,c,a) = \hat{q}^{c-r}\hat Q_a\in R.
\]

Now suppose $\la$ and $\nu$ are \bips of $n$ with $\la\doms\nu$. Let $G(\la,\nu)$ be the set of all pairs $(L,N)$ such that
\begin{itemize}
\item
$L$ is a rim hook of $\la$ and $N$ a rim hook of $\nu$,
\item
$[\la]\sm L=[\nu]\sm N$, and
\item
$\resi(\hand L)=\resi(\hand N)$.
\end{itemize}
Given a pair $(L,N)\in G(\la,\nu)$, define $\epsilon_{LN} = (-1)^{\lel(L)-\lel(N)}$, and let
\[
j_{\la\nu} = \prod_{(L,N)\in G(\la,\nu)}\big(\rez(\hand L)-\rez(\hand N)\big)^{\epsilon_{LN}}.
\]
Now for any pair of \bips $(\la,\mu)$ with $\mu$ restricted, we define
\[
J_{\la\mu} = \sum_{\nu\domsby\la}\nu_{\mathfrak{p}}(j_{\la\nu})\dn\nu\mu.
\]
The following statement is a special case of the \js formula \cite[Theorem 4.6]{jmjs}.

\begin{thm}[\xcite{jmjs}{Theorem 4.6}]
\label{jansch}
Suppose $\la$ and $\mu$ are \bips of $n$ with $\mu$ restricted.  Then the decomposition number $\dn\la\mu$ is at most $J_{\la\mu}$, and is non-zero if and only if $J_{\la\mu}$ is non-zero.
\end{thm}

Our particular choice of $\hat q,\hat Q_1,\hat Q_2$ means that the values $\nu_{\mathfrak{p}}(j_{\la\nu})$ will be easy to compute, and will not depend on the characteristic of $\bbf$ for the cases we will need to consider. In particular, suppose $(r,c,a)$ and $(s,d,b)$ are nodes. Then it is easy to check that:
\begin{itemize}
\item
if $a\neq b$, then
\[
\nu_{\fkp}(\rez(r,c,a)-\rez(s,d,b))=
\begin{cases}
1&\text{if }\resi(r,c,a)=\resi(s,d,b),
\\
0&\text{otherwise};
\end{cases}
\]
\item
if $a=b$, then
\[
\nu_{\fkp}(\rez(r,c,a)-\rez(s,d,b))=
\begin{cases}
1&\text{if }c-r=d-s\pm e,
\\
0&\text{if }c-r\nequiv d-s\ppmod e.
\end{cases}
\]

\end{itemize}

We can use the \js formula to refine \cref{basicdecomp,basicdecomp2}, by using a coarser order than the dominance order. We write $\la\blacktriangleright\nu$ if $\la\doms\nu$ and $\nu_{\fkp}(j_{\la\nu})\neq0$, and we extend $\blacktriangleright$ reflexively and transitively to give a partial order, which we call the \emph{\js dominance} order; note that this order depends on our fixed parameters $e$ and $\ka$. It is easy to see that the usual dominance order is a refinement of the \js dominance order, and that conjugation of \bips reverses the \js dominance order. Moreover, \cref{jansch} allows us to strengthen \cref{basicdecomp,basicdecomp2}.

\begin{propn}\label{basicdecompjs}
Suppose $\la$ and $\mu$ are \bips with $\mu$ restricted and $\dn\la\mu>0$, then $\mu\dm\dom\la\dom\mu$ in the \js dominance order.
\end{propn}

\begin{pf}
The fact that $\la\dom\mu$ is immediate from \cref{jansch}. This then implies that $\mu\dm\dom\la$, using \cref{conjdual} and the fact that the \js dominance order is reversed by conjugation.
\end{pf}

With this result in mind, we will use the \js order from now on, so we will write $\dom$ for the \js order. Although this order is harder to work with, it will help us to narrow down the cases we need to consider. Note that when trying to decide whether $\la\dom\mu$ in the \js order, the usual dominance order can be used as a ``first pass'': if $\la\ndom\mu$ in the usual dominance order, then certainly $\la\ndom\mu$ in the \js order. To examine the \js order more closely, it will be helpful to be able to visualise it on the abacus. Take two \bips $\la,\nu$, and take the \abds for these \bips. In order to obtain $\nu_{\fkp}(j_{\la\nu})\neq0$, we need to be able to get from the \abd for $\la$ to the \abd for $\nu$ by moving a bead from position $b_1$ to position $c_1$ on component $i_1$, and a bead from position $c_2$ to position $b_2$ on component $i_2$, where $i_1,i_2\in\{1,2\}$, $b_1-c_1=b_2-c_2>0$, and $b_1\equiv b_2\ppmod e$. We then have $\la\doms\nu$ \iff either $b_1>b_2$ or $i_1<i_2$. The leg length of the hook $L$ removed from $\la$ is the number of beads between positions $b_1$ and $c_1$, and the leg length of the hook removed from $\nu$ is the number of beads between positions $b_2$ and $c_2$.

For example, take $e=6$, $\ka=(\ol5,\ol4)$, $\la=\bp{2}{4,2,1^2}$ and $\nu=\bp{2,1^4}{4}$. Taking the bicharge $(11,10)$, we obtain the following \abds. (We draw the beads that have moved in white.)
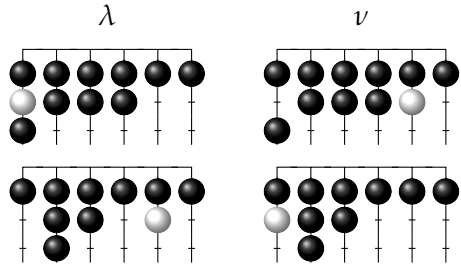
\[
\begin{array}{c@{\qquad}c}
\la&\nu
\\[5pt]
\abacus(lmmmmr,bbbbbb,obbbnn,bnnnnn)
&
\abacus(lmmmmr,bbbbbb,nbbbon,bnnnnn)
\\[18pt]
\abacus(lmmmmr,bbbbbb,nbbnon,nbnnnn)
&
\abacus(lmmmmr,bbbbbb,obbnnn,nbnnnn)
\end{array}
\]
Taking $b_1=b_2=10$, $c_1=c_2=6$, $i_1=1$ and $i_2=2$, we see that $j_{\la\nu}\neq0$ and $\la\domsby\nu$. Alternatively, we can see this by looking at the Young diagrams (in which we label nodes with their residues, and shade the involved rim hooks).
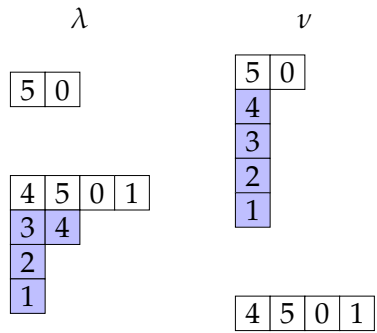
\[
\begin{array}{c@{\qquad}c}
\la&\nu
\\[5pt]
\begin{array}l
\gyoung(50,,,4501,!\blu34,2,1!\wht)
\end{array}
&
\begin{array}c
\gyoung(50,!\blu4,3,2,1,,,!\wht4501)
\end{array}
\end{array}
\]

\section{Blocks of weight $3$}\label{class3sec}

\subsection{Weight and the abacus}\label{wtabsec}

In this section we summarise results from \cite{mfweight} showing how to calculate the weight of a \bip from the abacus, and then use these to describe \bips of weight $3$.

Suppose $\la$ is a \bip, and construct the \abd for $\la$ using the bicharge $(k_1,k_2)$. We compute the weight of $\la$ recursively, in three stages.

First suppose there is a bead in the \abd with an unoccupied position immediately above. Define a new \abd (and hence a new \bip $\la^-$) by moving the bead up into this empty space (this corresponds to removing a rim $e$-hook from $\la$). Then \cite[Corollary~3.4]{mfweight} gives $\wt{\la^-}=\wt\la-2$.

Applying this repeatedly, we can assume that $\la$ is a \emph{bicore}; that is, it has no rim $e$-hooks. Now for each $x\in\zez$ define $\ga_x(\la)$ to equal the number of beads on runner $x$ in component $1$ minus the number of beads on runner $x$ in component $2$ of the \abd. (Note that $\ga_x(\la)$ depends on the choice of bicharge $(k_1,k_2)$, but the difference $\ga_x(\la)-\ga_y(\la)$ for $x,y\in\zez$ does not.) Now for $x,y\in\zez$ with $x\neq y$ define the bicore $s_{xy}(\la)$ by moving a bead from component $1$ to component $2$ on runner $x$, and a bead from component $2$ to component $1$ on runner $y$ of the \abd. Then \cite[Lemma 3.7]{mfweight} shows that $\wt{s_{xy}(\la)}=\wt\la-2(\ga_x(\la)-\ga_y(\la)-2)$.

This shows that whenever there are $x,y$ with $\ga_x(\la)-\ga_y(\la)\gs3$, we can replace $\la$ with a \bip $s_{xy}(\la)$ of smaller weight. By doing this repeatedly (and using the fact that weight is always non-negative) we can reduce to the case where $\la$ is a bicore with $\ga_x(\la)-\ga_y(\la)\ls2$ for all $x,y\in\zez$. In this case we use \cite[Proposition 3.8]{mfweight}, which gives $\wt\la=\min\{\card X,\card Y\}$, where
\begin{align*}
X&=\lset{x\in\zez}{\ga_x(\la)-\ga_y(\la)=2\text{ for some }y\in\zez},
\\
Y&=\lset{y\in\zez}{\ga_x(\la)-\ga_y(\la)=2\text{ for some }x\in\zez}.
\end{align*}
In particular, we note that $\wt\la=0$ \iff $\la$ is a bicore with $\ga_x(\la)-\ga_y(\la)\ls1$ for all $x,y\in\zez$.

\begin{eg}
Suppose $e=5$, $(\ka_1,\ka_2)=(\ol4,\ol4)$ and $(k_1,k_2)=(9,9)$, and let $\la$ be the \bip $\bp{5,3}{6,4,3}$, with the following \abd.
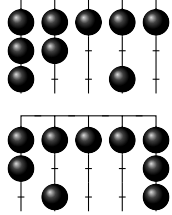
\[
\begin{array}{c}
\abacus(lmmmr,bbbbb,bbnnn,bnnbn)\\[18pt]
\abacus(lmmmr,bbbbb,bnnnb,nbnnb)
\end{array}
\]
Define a new \bip $\la^-=\bp{2,1}{6,2}$ by moving all beads as up as possible in their runners, and then define $\ol{\la}=s_{34}(\la^-)=\bp{2^2}{1^2}$.
\[
\begin{array}{c@{\qquad}c}
\la^-&\ol\la
\\[5pt]
\abacus(lmmmr,bbbbb,bbnbn,bnnnn)&\abacus(lmmmr,bbbbb,bbnnb,bnnnn)
\\[18pt]
\abacus(lmmmr,bbbbb,bbnnb,nnnnb)&\abacus(lmmmr,bbbbb,bbnbb,nnnnn)
\end{array}
\]
Then $\wt\la=\wt{\la^-}+4$ and $\wt{\la^-}=\wt{\ol{\la}}+2$, since $\ga_{\ol3}(\la^-)-\ga_{\ol4}(\la^-)=3$. Now $\ol{\la}$ satisfies $\ga_{\ol0}(\ol{\la})-\ga_{\ol3}(\ol{\la})=2$ and $\ga_x(\ol{\la})-\ga_y(\ol{\la})\ls1$ for all other $x,y\in\bbz/5\bbz$. In the above notation we have that $|X|=|Y|=1$ which implies $\wt{\ol{\la}}=1$. We conclude that $\wt\la=7$.
\end{eg}

We can use these results to describe blocks of weight $3$; this is very similar to the description of blocks of weight $2$ in \cite[Section 2.1]{mftypebwt2}. The above results show that if $\la$ is a \bip in a block $B$ of weight $3$, then one of the following occurs.
\begin{enumerate}
\item[1a.]
$\la$ has a rim $e$-hook, and removing this leaves a \bip of weight $1$.
\item[1b.]
$\la$ is a bicore, and there are $x,y\in\zez$ such that $\ga_x(\la)-\ga_y(\la)=3$ and $s_{xy}(\la)$ has weight $1$.
\item[2.]
$\la$ is a bicore with $\ga_i(\la)-\ga_j(\la)\ls2$ for all $i,j$, and $\min\{\card X,\card Y\}=3$, where $X$ and $Y$ are as defined above.
\end{enumerate}

We observe that if $\la^-$ is obtained by removing a rim $e$-hook from $\la$, then $\de_i(\la^-)=\de_i(\la)$ for all $i$, because $\de_i(\la)$ is determined by the total number of beads on each runner (across both components) in the \abd for $\la$, which does not change when we move a bead up its runner. Similarly, if $\la$ is a bicore then $\de_i(s_{xy}(\la))=\de_i(\la)$ for all $i$. In view of \cite[Lemma 3.3]{mfweight}, we then obtain $\card{s_{xy}(\la)}=n-e$. So in cases (1a) and (1b) there is a block $A$ of some $\hhh{n-e}$ with weight $1$ and with $\de_i(A)=\de_i(B)$ for all~$i$. 

If $B$ contains a \bip satisfying (2) above, then (following \cite{mfcore}) we say that $B$ is a \emph{core block}. By \cite[Theorem 3.1]{mfcore}, this is equivalent to the statement that there is no block $A$ of $\hhh m$ for any $m$ with weight less than $3$, and with $\de_i(A)=\de_i(B)$ for all $i$.

What this means is that if $B$ is a core block of weight $3$ then the \bips in $B$ all satisfy (2), while if $B$ is a non-core block then the \bips in $B$ all satisfy (1a) or (1b).

In \cref{inducsec} we will see that our main theorem is easy to prove for core blocks, but we need to consider non-core blocks in detail.

\subsection{The \bips in a non-core weight $3$ block}\label{bipsinblocksec}

Suppose $B$ is a non-core block of weight $3$. In this section we show how to construct all the \bips in $B$. Suppose $\la$ is a \bip in $B$, and take an \abd for $\la$ with bicharge $(k_1,k_2)$. Then we can construct a \bip $\nu$ lying in a block $A$ of weight $1$ with $\de_i(\nu)=\de_i(\la)$ for all $i$ (either by removing a rim $e$-hook or by applying $s_{xy}$ for some $x,y$). We call $A$ the \emph{weight $1$ block underlying $B$}. Now $\ga_x(\nu)-\ga_y(\nu)\ls2$ for all $x$ and $y$, and if we define the sets $X$ and $Y$ as in \cref{class3sec} (with $\nu$ in place of $\la$), then $\min\{\card X,\card Y\}=1$. We will assume that $\card Y=1$, with the other case being similar. We set $Z=X\cup Y$ and $\mathcal{C}(Z)=\zez\sm Z$; note in particular that $|Z|\gs2$. Now define a new \abd (with bicharge $(k_1+1,k_2-1)$) by moving a bead from component $2$ to component $1$ on runner $y$ (where $Y=\{y\}$).  We denote the resulting \bip $\xi_A$, and call it the \emph{\xr} of $A$. Then $\xi_A$ is independent of the choice of $\nu$ (and hence independent of the choice of $\la$), and the $(\ka_1+1,\ka_2-1)$-weight of $\xi_A$ is $0$. Moreover, if we define the integers $\ga_x(\xi_A)$ using the \abd with bicharge $(k_1+1,k_2-1)$, then
\[
\ga_x(\xi_A)-\ga_y(\xi_A)=
\begin{cases}
1&\text{if }x\in Z,\ y\in\mathcal{C}(Z)
\\
-1&\text{if }x\in\mathcal{C}(Z),\ y\in Z
\\
0&\text{otherwise}.
\end{cases}
\]
Now for each for each $z\in Z$ define a \bip $\nu_z$ in $A$ by taking the \abd for $\xi_A$ and moving a bead from component $1$ to component $2$ on runner $z$. Then the description of the \bips in a weight $1$ block in \cite[Section 4.2.1]{mfweight} shows that the \bips $\nu_z$ for $z\in X\cup Y$ are precisely the \bips in $A$. Now the results in \cref{class3sec} show that each \bip in $B$ is obtained from a \bip $\nu_z$ in $A$ either by adding a rim $e$-hook or by applying the operator $s_{xy}$ for some $x,y$ for which $\ga_x(\nu_z)-\ga_y(\nu_z)=1$. But the above formula for $\ga_x(\xi_A)-\ga_y(\xi_A)$ shows that $\ga_x(\nu_z)-\ga_y(\nu_z)=1$ \iff either $x\in Z\sm\{z\}$ and $y\in\mathcal{C}(Z)$, or $x\in\mathcal{C}(Z)$ and $y=z$.

This enables us to describe (and label) all the \bips in $B$. These come in three families, constructed from $\xi_A$ as follows.
\begin{itemize}
\item
For each $z\in Z$, $x\in\zez$ and $a\in\{1,2\}$, we obtain a \bip $\dl zxa$ from $\xi_A$ by moving a bead from component $1$ to component $2$ on runner $z$, then moving the lowest bead on runner $x$ on component $a$ down one position.
\item
For each $x\in\mathcal{C}(Z)$, we obtain a \bip $\dw x$ from $\xi_A$ by moving a bead from component $1$ to component $2$ on runner $x$.
\item
For each $w,z\in Z$ with $w\neq z$ and for each $y\in\mathcal{C}(Z)$, we obtain a \bip $\ddu wzy$ from $\xi_A$ by moving beads from component $1$ to component $2$ on runners $w$ and $z$, and then moving a bead from component $2$ to component $1$ on runner $y$.
\end{itemize}

\begin{eg}
Suppose $e=4$, $(\ka_1,\ka_2)=(\ol0,\ol3)$ and $(k_1,k_2)=(8,7)$. Let $B$ be the weight-$3$ block of $\hhh{10}$ containing the \bip $\la=\bp{4}{4,1^2}$. Then the underlying weight-$1$ block $A$ comprises the \bips $\bp{3,1^3}{\vn}$, $\bp{3}{1^3}$ and $\bp{\vn}{4,1^2}$, and the \xr of $A$ is $\xi_A=\bp{2}{1^2}$. These \bips have the following \abds.
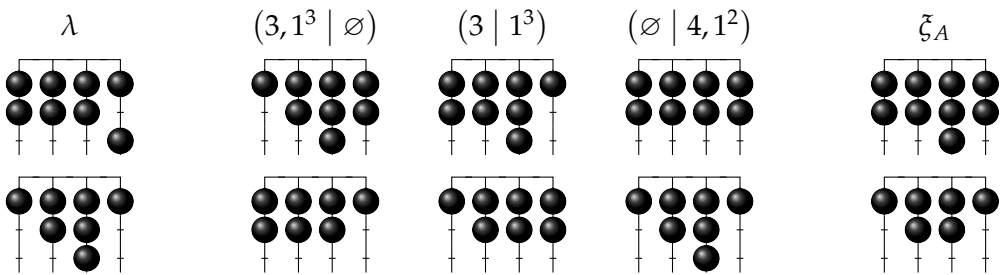
\[
\begin{array}{c@{\qquad\qquad}c@{\qquad}c@{\qquad}c@{\qquad\qquad}c}
\la&\bp{3,1^3}{\vn}&\bp{3}{1^3}&\bp\vn{4,1^2}&\xi_A
\\[5pt]
\abacus(lmmr,bbbb,bbbn,nnnb)&\abacus(lmmr,bbbb,nbbb,nnbn)&\abacus(lmmr,bbbb,bbbn,nnbn)&\abacus(lmmr,bbbb,bbbb,nnnn)&\abacus(lmmr,bbbb,bbbb,nnbn)
\\[18pt]
\abacus(lmmr,bbbb,nbbn,nnbn)&\abacus(lmmr,bbbb,bbbn,nnnn)&\abacus(lmmr,bbbb,nbbb,nnnn)&\abacus(lmmr,bbbb,nbbn,nnbn)&\abacus(lmmr,bbbb,nbbn,nnnn)
\end{array}
\]
We see that $Z=\{\ol0,\ol2,\ol3\}$. So there are $28$ \bips in $B$: the \bips $\dl zxa$ for $z\in\{\ol0,\ol2,\ol3\}$, $x\in\bbz/4\bbz$ and $a\in\{1,2\}$, together with the four \bips $\dw1$, $\ddu021$, $\ddu031$ and $\ddu231$. For example, the \bip $\la$ above is $\dl231$, while $\ddu031$ is the \bip $\bp{3^2,1^2}{1^2}$ shown below.
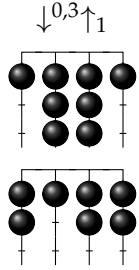
\[
\begin{array}{c}
\ddu031
\\[5pt]
\abacus(lmmr,bbbb,nbbn,nbbn)
\\[18pt]
\abacus(lmmr,bbbb,bnbb,nnnn)
\end{array}
\]
\end{eg}

\subsection{Inductive approach to the main theorem}\label{inducsec}

Now we give some results to facilitate a proof of \cref{main} by induction. We begin with an analogue of a result for partitions in \cite{jmsmalld}, which greatly restricts the set of \bips we need to consider. We give the result for blocks of weight $3$ only, but a more general result along the lines of \cite[Proposition~2.1]{jmsmalld} can easily be proved.

\begin{lemma}\label{simpledn0}
Suppose $B$ is a weight $3$ block of $\hhh n$, and that \cref{main} holds with $n$ replaced by any smaller integer. Suppose $\la,\mu$ are \bips of $n$ with $\mu$ restricted, and $i\in\zez$.
\begin{enumerate}
\item
If $\nor i\mu\gs\max\{1,\rem i\la\}$, then $\dn\la\mu\ls1$.
\item
If $\nor i\mu\gs\max\{1,\de_i(B)+1\}$, then $\dn\la\mu\ls1$.
\end{enumerate}
\end{lemma}

\begin{pf}
Recall the notation $\ep_iM$ from \cref{branchsec}, and that (from \cref{brchrules}) $\ep_i\spe\la=\rem i\la$ and $\ep_i\smp\mu=\nor i\mu$.
\begin{enumerate}
\item
Since $\ee i$ is an exact functor, if $M$ is a module and $N$ is a subquotient of $M$, then $\ep_iN\ls\ep_iM$. In particular, if $\nor i\mu>\rem i\la$, then $\ep_i\smp\mu>\ep_i\spe\la$, and therefore $\smp\mu$ cannot be a composition factor of $\spe\la$.

Suppose instead that $\nor i\mu=\rem i\la=r$, say. Then $r=\ep_i\spe\la=\ep_i\smp\mu$. Now $\ee i^{(r)}(\spe\la)\cong\spe{\la^-}$, and $\ee i^{(r)}(\smp\mu)\cong\smp{\mu^-}$, where $\la^-$ is obtained by removing all the removable $i$-nodes from $\la$, and $\mu^-$ is obtained by removing all the normal $i$-nodes from $\mu$. Therefore
\[
\dn{\la^-}{\mu^-}=[\ee i^{(r)}\spe\la:\smp{\mu^-}]=\sum_\nu\dn\la\nu[\ee i^{(r)}\smp\nu:\smp{\mu^-}]\gs\dn\la\mu,
\]
summing over restricted \bips $\nu$ of $n$. Because $r=\rem i\la$, we obtain $r\gs\de_i(\la)$. So by \cref{36lem} the weight of $\la^-$ is at most $3$. If $\la^-$ has weight $3$, then by the inductive assumption $\dn{\la^-}{\mu^-}\ls1$. If $\la^-$ has weight less than $3$, then $\dn{\la^-}{\mu^-}\ls1$ from \cref{w012}. Either way, we obtain $\dn\la\mu\ls1$.
\item
We show that the number $r=\rem i\la$ is at most $\max\{1,\de_i(B)+1\}$; then in particular $r\ls\nor i\mu$, and we can use part (1) of the \lcnamecref{simpledn0}. Suppose for a contradiction that $r\gs2$ and $r\gs\de_i(B)+2$, and let $\la^-$ be obtained by removing all the removable $i$-nodes from $\la$. Then by \cref{36lem} the weight of $\la^-$ equals $3-r(r-\de_i(B))$. Since by assumption $r\gs2$ and $r-\de_i(B)\gs2$, this makes the weight of $\la^-$ negative, a contradiction.\qedhere
\end{enumerate}
\end{pf}

We deduce two very useful corollaries.

\begin{cory}\label{casei}
Suppose $B$ is a weight $3$ block of $\hhh n$ with $\de_i(B)\ls0$ for all $i\in\zez$, and that \cref{main} holds with $n$ replaced by any smaller integer. Then \cref{main} holds for $B$.
\end{cory}

\begin{pf}
Suppose $\la,\mu$ are \bips in $B$ with $\mu$ restricted. Then $\smp\mu\downarrow_{\hhh{n-1}}\neq0$, so $\ee i\smp\mu\neq0$ for some $i\in\zez$. This means that $\nor i\mu\gs1$, so $\nor i\mu\gs\max\{1,\de_i(B)+1\}$, and therefore $\dn\la\mu\ls1$ by \cref{simpledn0}(2).
\end{pf}

\begin{cory}\label{untricky}
Suppose $B$ is a weight $3$ block of $\hhh n$, and that \cref{main} holds with $n$ replaced by any smaller integer. If $\la$ is a \bip in $B$ and $i\in\zez$ such that $\de_i(B)\gs1$ and $\la$ has no addable $i$-nodes, then $\dn\la\mu\ls1$ for all $\mu$.
\end{cory}

\begin{pf}
The assumption that $\la$ has no addable $i$-nodes means that $\rem i\la=\de_i(B)$, while $\nor i\mu\gs\de_i(B)$. So $\nor i\mu\gs\max\{1,\rem i\la\}$, and the result follows from \cref{simpledn0}(1).
\end{pf}

So to prove \cref{main} by induction using \cref{untricky}, it suffices to consider only those $\la$ that have at least one addable $i$-node for \emph{every} $i$ for which $\de_i(B)\gs1$. Call such $\la$ \emph{\tricky}. The next few lemmas will show that in most cases $B$ has no \tricky \bips.

\begin{lemma}\label{trickycore}
Suppose $B$ is a core block, and that $\de_i(B)\gs1$ for some $i\in\zez$. Then $B$ contains no \tricky \bips.
\end{lemma}

\begin{pf}
Suppose $\la$ is a \bip in $B$. Then we claim that $\la$ has no addable $i$-nodes. Assume for a contradiction that $\la$ has at least one addable $i$-node. Then because $\de_i(B)\gs1$, we know that $\la$ has at least two removable $i$-nodes. The two components of $\la$ are both $e$-cores, and an $e$-core cannot have addable and removable nodes of the same residue, so the addable $i$-nodes of $\la$ are in one component and the removable $i$-nodes in the other. We assume that the removable $i$-nodes are in component $1$ (the other case is similar). This means that in an \abd for $\la$ there are at least two positions $b$ on runner $i$ such that there is a bead in position $b$ but not in position $b-1$. So there are at least two more beads on runner $i$ than on runner $i-1$ in component $1$ (or at least three more, in the case $i=\ol0$). Similarly, in component $2$ of the \abd, there is at least one position $b$ on runner $i-1$ such that there is a bead at position $b$ but no bead at position $b+1$. So there are more beads on runner $i-1$ than on runner $i$ in the component $2$ (or at least as many beads on runner $i-1$ as on runner $i$, in the case $i=\ol0$). We deduce that $\ga_i(\la)-\ga_{i-1}(\la)\gs3$. But the assumption that $B$ is a core block gives $\ga_i(\la)-\ga_{i-1}(\la)\ls2$, a contradiction. 
\end{pf}

In view of \cref{trickycore}, it suffices to consider non-core blocks from now on.

\begin{lemma}\label{k3}
Suppose $\de_i(B)\gs3$ for some $i\in\zez$. Then $B$ contains no \tricky \bips.
\end{lemma}

\begin{pf}
If $\la$ is a \bip in $B$, then we claim that $\la$ cannot have an addable $i$-node: if it did, then by \cref{36lem} the weight of a \bip obtained by adding an addable $i$-node to $\la$ would be $3-(\de_i(\la)+1)$, which is negative, a contradiction.
\end{pf}

\begin{lemma}\label{kij}
Suppose there are $i,j\in\zez$ such that $\de_i(B),\de_j(B)\gs1$, and $j\neq i,i\pm1$. Then $B$ contains no \tricky \bips.
\end{lemma}

\begin{pf}
Suppose $\la$ is an \tricky \bip in $B$. Then $\la$ has both an addable $i$-node $\fkm$ and an addable $j$-node $\fkn$. Define \bips $\mu,\nu$ by $[\mu]=[\la]\cup\{\fkm\}$ and $[\nu]=[\mu]\cup\{\fkn\}$. \cref{36lem} then gives $\wt\mu=3-(\de_i(B)+1)\ls1$. Because $j\neq i,i\pm1$, adding the $i$-node to $\la$ does not affect the addable or removable $j$-nodes, so $\de_j(\mu)=\de_j(\la)$. Hence $\wt\nu=\wt\mu-(\de_j(B)+1)\ls-1$, a contradiction.
\end{pf}

\begin{lemma}\label{consec12}
Suppose there are $i,j\in\zez$ such that $i\neq j$,  $\de_i(B)\gs1$ and $\de_j(B)\gs2$. Then $B$ contains no \tricky \bips.
\end{lemma}

\begin{pf}
If $j\neq i\pm1$ then the result follows from \cref{kij}, so assume that $j=i+1$ (the case $j=i-1$ is similar). Suppose $\la$ is an \tricky \bip in $B$. Then $\la$ has both an addable $i$-node $\fkm$ and an addable $(i+1)$-node $\fkn$. Define \bips $\mu,\nu$ by $[\mu]=[\la]\cup\{\fkm\}$ and $[\nu]=[\mu]\cup\{\fkn\}$. \cref{36lem} now gives $\wt\mu=3-(\de_i(B)+1)\ls1$. Now from \cite[Lemma 3.1]{mfweight} we obtain $\de_{i+1}(\mu)=\de_{i+1}(\la)-1$, so that $\wt\nu=\wt\mu-\de_{i+1}(B)\ls-1$, and again we have a contradiction.
\end{pf}

These lemmas mean that there are only four types of non-core block $B$ in which there can be \tricky \bips. These are exactly analogous to the four types of weight $3$ blocks of symmetric groups considered in \cite[Section 4]{jmsmalld}. These types can be defined according to the values $\de_i(B)$, as follows.

\begin{enumerate}[label=\Roman*.]
\item
$\de_i(B)\ls0$ for all $i\in\zez$.
\item
There is $i\in\zez$ such that $\de_i(B)=1$ while $\de_j(B)\ls0$ for all $j\neq i$.
\item
There is $i\in\zez$ such that $\de_i(B)=\de_{i+1}(B)=1$ while $\de_j(B)\ls0$ for all $j\neq i,i+1$.
\item
There is $i\in\zez$ such that $\de_i(B)=2$ while $\de_j(B)\ls0$ for all $j\neq i$.
\end{enumerate}

We know from \cref{casei} that \cref{main} holds for blocks of type I. To prove \cref{main}, we need to show that if $B$ is a block of $\hhh n$ of type II, III or IV and if \cref{main} holds with $n$ replaced by any smaller integer, then it holds for $B$. We address the three types of block in \cref{type2sec,type3sec,type4sec}.

\section{Blocks of type II}\label{type2sec}

Throughout this section, we assume $B$ is a block of $\hhh n$ and $i\in\bbz$ is such that $\de_i(B)=1$ while $\de_j(B)\ls0$ for all $j\nequiv i\ppmod e$.  Our aim is to show that if \cref{main} holds with $n$ replaced by any smaller integer, then it holds for $B$.

\subsection{The abacus for a block of type II}

Let $\xi$ be the \xr of $B$. Whenever we refer to the \abd for $\xi$, we will mean the $(k_1+1,k_2-1)$-\abd, and whenever we refer to the residues of nodes for $\xi$, we mean with respect to $(\kappa_1+1,\kappa_2-1)$. With this convention in mind, we have $\de_i(\xi)=1$ while $\de_j(\xi)\ls0$ for $j\nequiv i\ppmod e$. This means that $\xi$ has exactly one removable node, which has residue $i$ (because by \cref{36lem} a \bip of weight $0$ cannot have an addable and a removable node of the same residue). We will assume that this removable node lies in component $1$ (in the opposite case we can replace $\ka$ with $(\ka_2,\ka_1)$ and appeal to the results in \cref{conjswsec}).

By considering the \abd for $\xi$ and using the fact that $\xi$ has $(\ka_1+1,\ka_2-1)$-weight $0$, we find that there are integers $i\ls j\ls k\ls l\ls e+i-2$ such that in component $1$ of $\xi$ there is a removable $i$-node, and addable nodes of residues $j+1$ and $l+1$, while in component $2$ there is just an addable node of residue $k+1$ (so in fact $\ka_1=\widebar{j+l+1-i}$, and $\ka_2=\widebar{k+2}$). We therefore have $Z=\{\ol i,\dots,\ol j\}\cup\{\widebar{k+1},\dots,\ol l\}$. The fact that $|Z|\gs2$ implies that either $i<j$ or $k<l$.

\begin{eg}
Suppose $e=11$, $\kappa=(\ol0,\ol7)$ and $\xi=\bp{3^3}\vn$, with the following abacus display.
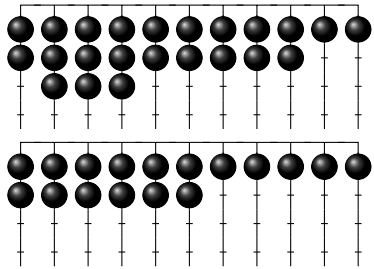
\[
\begin{array}c
\abacus(lmmmmmmmmmr,bbbbbbbbbbb,bbbbbbbbbnn,nbbbnnnnnnn,nnnnnnnnnnn)
\\[20pt]
\abacus(lmmmmmmmmmr,bbbbbbbbbbb,bbbbbbnnnnn,nnnnnnnnnnn,nnnnnnnnnnn)
\end{array}
\]
We can take $(i,j,k,l)=(1,3,5,8)$, and $Z=\{\ol1,\ol2,\ol3,\ol6,\ol7,\ol8\}$.
\end{eg}

\subsection{\Tricky \bips in $B$}

Suppose $\la$ is a \bip in $B$. From the results in \cref{inducsec}, we need only consider Specht modules labelled by \tricky \bips; that is, those with addable $i$-nodes. By examining the \abds of the \bips in $B$, we can easily classify these. There are $3\card Z$ \tricky \bips, which we label $\al_x,\be_x,\ga_x$ for $x\in Z$ as follows:
\[
\al_x=
\begin{cases}
\dl xi2&\text{ if }x\neq i
\\
\dw{i-1}&\text{ if }x=i,
\end{cases}
\qquad
\be_x=
\begin{cases}
\dl x{i-1}2&\text{ if }x\neq i
\\
\dl ii1&\text{ if }x=i,
\end{cases}
\qquad
\ga_x=
\begin{cases}
\ddu ix{i-1}&\text{ if }x\neq i
\\
\dl i{i-1}1&\text{ if }x=i.
\end{cases}
\]

It will be helpful to write these \bips explicitly:
{
\medmuskip=2mu
\begin{align*}
\al_x&=
\begin{cases}
\bp{(e+i-l)^{j-i+1},1^{l-i+1}}{e+i-k-2}&\text{if }x=i
\\
\hbox to 350pt{$\bp{(e+i-l)^{j-x},(e+i-l-1)^{x-i}}{e+x-k-1,e+i-k,1^{k-i}}$\hfil}&\text{if }i<x\ls j
\\
\bp{(e+i-l)^{j-i+1},1^{l-x}}{e+i-k-1,x-k,1^{k-i}}&\text{if }k<x\ls l,
\end{cases}
\\
\be_x&=
\begin{cases}
\bp{(e+i-l)^{j-i+1},1^{l-i}}{e+i-k-1}&\text{if }x=i
\\
\hbox to 350pt{$\bp{(e+i-l)^{j-x},(e+i-l-1)^{x-i}}{e+x-k-1,e+i-k-1,1^{k-i+1}}$\hfil}&\text{if }i<x\ls j
\\
\bp{(e+i-l)^{j-i+1},1^{l-x}}{e+i-k-2,x-k,1^{k-i+1}}&\text{if }k<x\ls l,
\end{cases}
\\
\ga_x&=
\begin{cases}
\bp{(e+i-l)^{j-i},e+i-l-1,1^{l-i+1}}{e+i-k-1}&\text{if }x=i
\\
\hbox to 350pt{$\bp{(e+i-l)^{j-x},(e+i-l-1)^{x-i-1},e+i-l-2}{e+x-k-1,e+i-k,1^{k-i+1}}$\hfil}&\text{if }i<x\ls j
\\
\bp{(e+i-l)^{j-i},e+i-l-1,1^{l-x}}{e+i-k-1,x-k,1^{k-i+1}}&\text{if }k<x\ls l.
\end{cases}
\end{align*}
}
Each \tricky \bip has exactly one addable $i$-node, and the labelling is chosen so that $\al_x,\be_x,\ga_x$ all yield the same \bip when the addable $i$-node is added, with $\al_x\doms\be_x\doms\ga_x$.

\subsection{Proving \cref{main} in $B$}\label{proving2}

Our aim is to show that $\dn\la\mu\ls1$ for all \bips $\la,\mu$ in $B$, given the inductive assumption. We will do this by looking at one restricted \bip $\mu$ at a time. The next lemma will help in narrowing down the cases we need to consider.

\begin{lemma}\label{domcond}
Suppose $\mu$ is a restricted \bip in $B$, and that $\dn\la\mu>1$ for some $\la$. Then $\nor i\mu=1$, and there is some $x\in Z$ for which $\dn{\al_x}\mu$, $\dn{\be_x}\mu$, $\dn{\ga_x}\mu$ are all non-zero; in particular, $\mu\dm\dom\al_x,\be_x,\ga_x\dom\mu$.
\end{lemma}

\begin{pf}
If $\la$ is such that $\dn\la\mu>1$, then $\la$ is \tricky, so $\la\in\{\al_x,\be_x,\ga_x\}$ for some $x$. The conclusion that $\mu$ has only one normal node comes from \cref{simpledn0}(2).

Now we claim that $\dn{\al_x}\mu,\dn{\be_x}\mu,\dn{\ga_x}\mu\gs1$, which gives the result using \cref{basicdecompjs}. Each of $\al_x,\be_x,\ga_x$ has two removable $i$-nodes, and the branching rules show that there are \bips $\hat\al_x,\hat\be_x,\hat\ga_x$ of $n-1$ such that
\begin{align*}
\ee i\spe{\al_x}&\sim\spe{\hat\al_x}+\spe{\hat\be_x},
\\
\ee i\spe{\be_x}&\sim\spe{\hat\al_x}+\spe{\hat\ga_x},
\\
\ee i\spe{\ga_x}&\sim\spe{\hat\be_x}+\spe{\hat\ga_x}.
\end{align*}
Because $\nor i\mu=1$, there is a restricted \bip $\tilde\mu$ such that $\ee i\smp\mu\cong\smp{\tilde\mu}$. Now (analogously to the proof of \cite[Proposition 2.8(2)]{mfwt3}) we obtain
\[
\dn{\al_x}\mu+\dn{\hat\ga_x}{\tilde\mu}=\dn{\be_x}\mu+\dn{\hat\be_x}{\tilde\mu}=\dn{\ga_x}\mu+\dn{\hat\al_x}{\tilde\mu}.
\]
Since by assumption the decomposition numbers for weight $3$ blocks of $\hhh{n-1}$ are all $0$ or $1$, this means that the decomposition numbers $\dn{\al_x}\mu$, $\dn{\be_x}\mu$, $\dn{\ga_x}\mu$ differ by at most $1$. So if one of them is greater than $1$, then they are all positive.
\end{pf}

\cref{domcond} greatly restricts the set of restricted \bips $\mu$ that we need to consider. To apply \cref{domcond} we look at the \js dominance order among the \tricky \bips, which is easy to check given the explicit description above: if $\epsilon$ is any of the symbols $\al,\be,\ga$, then
\[
\epsilon_i\doms\epsilon_{k+1}\doms\epsilon_{k+2}\doms\cdots\doms\epsilon_l\doms\epsilon_{i+1}\doms\epsilon_{i+2}\doms\cdots\doms\epsilon_j.
\]
In particular, the most dominant $\ga_x$ is $\ga_i$. So to show that $\dn\la\mu\ls1$ for all $\mu$, we need only consider \bips $\mu$ for which:
\begin{enumerate}
\item\label{muklesh}
$\mu$ is restricted;
\item\label{gadom}
$\ga_i\dom\mu$; and
\item\label{normi}
$\nor i\mu=1$, and $\nor j\mu=0$ for all $j\neq i$.
\end{enumerate}

Analysis of the possible abacus configurations enables a classification of the \bips $\mu$ satisfying (\ref{gadom}) and (\ref{normi}). By working through the different \bips in $B$, we can check that the only \bips $\mu$ satisfying these conditions are $\dl i{i-1}2$ and $\dl ii2$.

(Note that if $i<k$ then the \bip $\mu=\dl ik2$ does not satisfy (\ref{gadom}): in this case
\begin{align*}
\mu&=\bp{(e+i-l)^{j-i}}{e-1,e+i-k},
\\
\ga_i&=\bp{(e+i-l)^{j-i},e+i-l-1,1^{l-i+1}}{e+i-k-1},
\end{align*}
so $\ga_i\doms\mu$ in the usual dominance order. But by examining residues we can check that $\ga_i\ndoms\mu$ in the \js dominance order.)

Now we introduce condition (1). First, from \cref{resterest} we know that $\dl ii2$ is not restricted, since its second component is not $e$-restricted. So it remains to check when $\mu=\dl i{i-1}2$ is restricted. We do this by removing good nodes, using the results of \cref{kleshbipsec}.

We begin by drawing an \abd for $\mu$ (where we adjust the numbers of beads in order to make runner $i-1$ the leftmost runner).
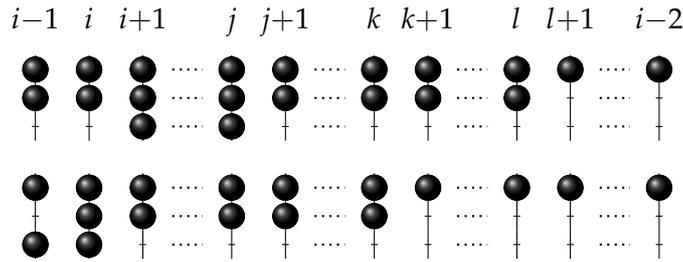
\[
\begin{array}{ccc@{\ \ }c@{\ \ }cc@{\ \ }c@{\ \ }cc@{\ \ }c@{\ \ }cc@{\ \ }c@{\ \ }c}
\mathclap{i{-}1}&\mathclap{i}&\mathclap{i{+}1}&&\mathclap{j}&\mathclap{j{+}1}&&\mathclap{k}&\mathclap{k{+}1}&&\mathclap{l}&\mathclap{l{+}1}&&\mathclap{i{-}2}
\\[6pt]
\abacus(b,b,n)&
\abacus(b,b,n)&
\abacus(b,b,b)&
\abacus(h,h,h)&
\abacus(b,b,b)&
\abacus(b,b,n)&
\abacus(h,h,h)&
\abacus(b,b,n)&
\abacus(b,b,n)&
\abacus(h,h,h)&
\abacus(b,b,n)&
\abacus(b,n,n)&
\abacus(h,h,h)&
\abacus(b,n,n)
\\[20pt]
\abacus(b,n,b)&
\abacus(b,b,b)&
\abacus(b,b,n)&
\abacus(h,h,h)&
\abacus(b,b,n)&
\abacus(b,b,n)&
\abacus(h,h,h)&
\abacus(b,b,n)&
\abacus(b,n,n)&
\abacus(h,h,h)&
\abacus(b,n,n)&
\abacus(b,n,n)&
\abacus(h,h,h)&
\abacus(b,n,n)
\end{array}
\]
Now we remove good nodes of residues $i,i+1,\dots,k,\ i-1,i-2,\dots,k+1$ in turn, to obtain the \bip $\xi=\bp{(e+i-l)^{j-i}}{e+i-k-1}$.

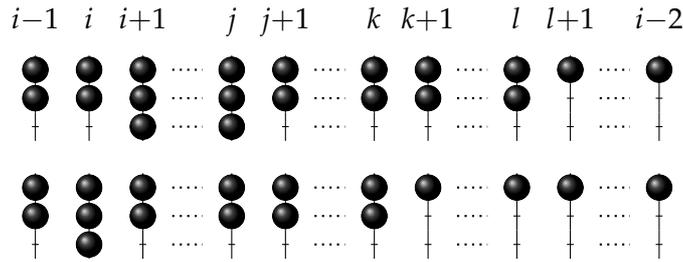
\[
\begin{array}{ccc@{\ \ }c@{\ \ }cc@{\ \ }c@{\ \ }cc@{\ \ }c@{\ \ }cc@{\ \ }c@{\ \ }c}
\mathclap{i{-}1}&\mathclap{i}&\mathclap{i{+}1}&&\mathclap{j}&\mathclap{j{+}1}&&\mathclap{k}&\mathclap{k{+}1}&&\mathclap{l}&\mathclap{l{+}1}&&\mathclap{i{-}2}
\\[6pt]
\abacus(b,b,n)&
\abacus(b,b,n)&
\abacus(b,b,b)&
\abacus(h,h,h)&
\abacus(b,b,b)&
\abacus(b,b,n)&
\abacus(h,h,h)&
\abacus(b,b,n)&
\abacus(b,b,n)&
\abacus(h,h,h)&
\abacus(b,b,n)&
\abacus(b,n,n)&
\abacus(h,h,h)&
\abacus(b,n,n)
\\[20pt]
\abacus(b,b,n)&
\abacus(b,b,b)&
\abacus(b,b,n)&
\abacus(h,h,h)&
\abacus(b,b,n)&
\abacus(b,b,n)&
\abacus(h,h,h)&
\abacus(b,b,n)&
\abacus(b,n,n)&
\abacus(h,h,h)&
\abacus(b,n,n)&
\abacus(b,n,n)&
\abacus(h,h,h)&
\abacus(b,n,n)
\end{array}
\]
This \bip lies in the weight $1$ block underlying $B$, and (from the results in \cref{kleshbipsec}) is restricted \iff $\mu$ is. Hence we need to understand whether $\xi$ is restricted or not. \cite[Lemma 4.9]{mfweight} shows that a \bip of weight $1$ is restricted \iff it is not the most dominant \bip in its block. Using the explicit description of weight 1 blocks \cite[Theorem 4.4]{mfweight}, we find all the \bips in the same block as $\xi$ and we see that $\xi$ is the most dominant \bip in its block \iff $k=l$. It follows that $\mu$ is restricted \iff $k<l$.

We deal with this case using the \js formula. From \cref{domcond} we need only consider \tricky \bips $\al_x,\be_x,\ga_x$ for which $\ga_x\dom\mu$, and by looking at the explicit expressions for the \tricky \bips, we find that this happens \iff $x=i$ or $k<x\ls l$. So to show that $\dn\la\mu\ls1$ for all $\la$ it suffices to show that
\[
\dn{\ga_i}\mu=\dn{\ga_{k+1}}\mu=\dots=\dn{\ga_{l-1}}\mu=\dn{\be_l}\mu=0.
\]
To do this using the \js order, we first find all the \bips $\la$ such that $\be_l\dom\la\dom\mu$ or $\ga_x\dom\la\dom\mu$ for some $x\in\{i\}\cup\{k+1,\dots,l\}$. In fact by examining the explicit forms of these \bips, we easily find that the only such \bips $\la$ are $\mu$, $\be_l$, $\ga_i$ and $\ga_{k+1},\dots,\ga_l$, with the \js order on these \bips given by the following diagram.
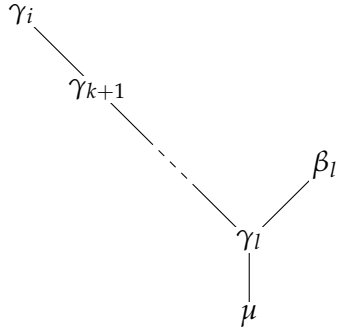
\[
\begin{tikzpicture}[scale=1,every node/.style={fill=white,inner sep=1pt}]
\draw(0,0)coordinate(mu);
\draw(0,1)coordinate(gal);
\draw(1,2)coordinate(bel);
\draw(-.67,1.67)coordinate(gaa);
\draw(-1.33,2.33)coordinate(gab);
\draw(-2,3)coordinate(gak);
\draw(-3,4)coordinate(gai);
\draw(mu)--(gal)--(bel);
\draw(gal)--(gaa);
\draw[dashed](gaa)--(gab);
\draw(gab)--(gai);
\draw(mu)node{$\mu$};
\draw(gal)node{$\ga_l$};
\draw(bel)node{$\be_l$};
\draw(gak)node{$\ga_{k+1}$};
\draw(gai)node{$\ga_i$};
\end{tikzpicture}
\]
The \js coefficients we need are given by
\[
j_{\mu\ga_x}=
\begin{cases}
(-1)^{l-x}&(k<x\ls l)
\\
(-1)^{l-k}&(x=i),
\end{cases}
\qquad
j_{\mu\be_l}=-1,\qquad j_{\ga_l\ga_x}=
\begin{cases}
(-1)^{l-x+1}&(k<x<l)
\\
(-1)^{l-k+1}&(x=i),
\end{cases}
\qquad
j_{\ga_l\be_l}=1.
\]
Now we deduce that $\dn{\ga_l}\mu=1$, and then (successively) $\dn{\ga_{l-1}}\mu=\dots=\dn{\ga_{k+1}}\mu=\dn{\ga_i}\mu=0$ as well as $\dn{\be_l}\mu=0$.

\section{Blocks of type III}\label{type3sec}

Throughout this section, we assume $B$ is a block of $\hhh n$ and $i\in\bbz$ is such that $\de_i(B)=\de_{i+1}(B)=1$ while $\de_j(B)\ls0$ for all $j\nequiv i,i+1\ppmod e$. Our aim is to show that if \cref{main} holds with $n$ replaced by any smaller integer, then it holds for $B$.

\subsection{The abacus for a block of type III}

Let $\xi$ be the \xr of $B$. Whenever we refer to the \abd for $\xi$, we will mean the $(k_1+1,k_2-1)$-\abd, and whenever we refer to the residues of nodes for $\xi$, we mean with respect to $(\kappa_1+1,\kappa_2-1)$. With this convention in mind, we have $\de_i(\xi)=\de_{i+1}(\xi)=1$ while $\de_j(\xi)\ls0$ for $j\nequiv i,i+1\ppmod e$. This means that $\xi$ has exactly two removable nodes, of residues $i$ and $i+1$. We will assume that the removable $(i+1)$-node lies in component $1$ (in the opposite case we can replace $\ka$ with $(\ka_2,\ka_1)$ and appeal to the results in \cref{conjswsec}). This means in particular that $\ga_{i+1}(\xi)=\ga_i(\xi)+1$. Now we claim that the removable $i$-node of $\xi$ must lie in component $2$; if not, then we obtain $\ga_i(\xi)=\ga_{i-1}(\xi)+1$ and hence $\ga_{i+1}(\xi)-\ga_{i-1}(\xi)=2$, and now the results in \cref{wtabsec} show that $\xi$ has positive weight, a contradiction.

Now we can describe the \abd for $\xi$ more precisely: there are integers $i+1\ls j\ls k\ls l\ls m\ls e+i-2$ such that in component $1$ of $\xi$ there is a removable node of residue $i+1$ and addable nodes of residues $j+1$ and $l+1$, while in component $2$ there is a removable node of residue $i$ and addable nodes of residues $k+1$ and $m+1$. (So in fact $\kappa_1=\widebar{j+l-i}$, and $\kappa_2=\widebar{k+m+3-i}$). We therefore have $Z=\{\ol {i+1},\dots,\ol j\}\cup\{\widebar {k+1},\dots,\widebar l\}\cup\{\widebar {m+1},\dots,\ol {i-1}\}$.

\begin{eg}
Suppose $e=14$, $\kappa=(\ol{10},\ol5)$ and $\xi=\bp{7^2}{3^6}$, with the following abacus display.
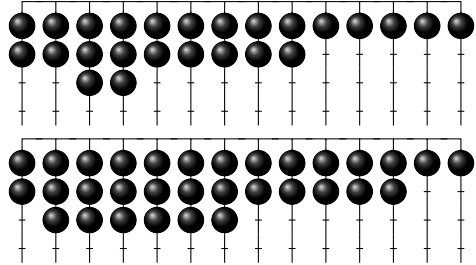
\[
\begin{array}c
\abacus(lmmmmmmmmmmmmr,bbbbbbbbbbbbbb,bbbbbbbbbnnnnn,nnbbnnnnnnnnnn,nnnnnnnnnnnnnn)
\\[20pt]
\abacus(lmmmmmmmmmmmmr,bbbbbbbbbbbbbb,bbbbbbbbbbbbnn,nbbbbbbnnnnnnn,nnnnnnnnnnnnnn)
\end{array}
\]
We can take $(i,j,k,l,m)=(1,3,6,8,11)$, and $Z=\{\ol0,\ol2,\ol3,\ol7,\ol8,\widebar{12},\widebar{13}\}$.
\end{eg}

\subsection{\Tricky \bips in $B$}

The \tricky \bips in $B$ are those with both addable $i$-nodes and addable $(i+1)$-nodes. By examining the \abds of all the \bips in $B$, we find that there are four \tricky \bips, which we label as follows, showing the configuration of runners $i-1,i,i+1$ of their \abds:
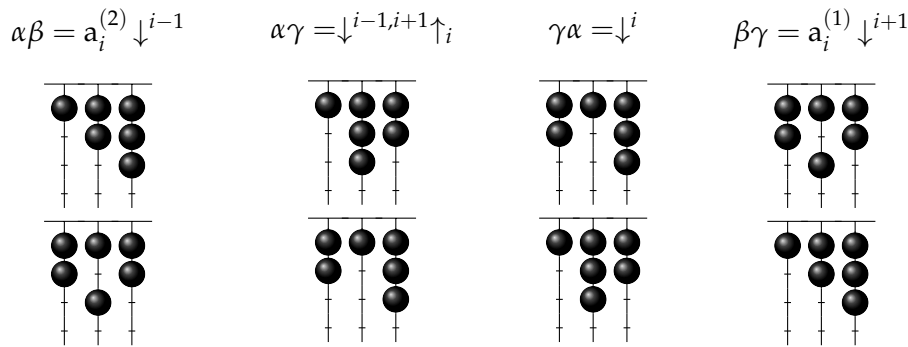
\[
\begin{array}c
\albe=\dl{i-1}i2
\\[12pt]
\abacus(mmm,bbb,nbb,nnb,nnn)
\\[20pt]
\abacus(mmm,bbb,bnb,nbn,nnn)
\end{array}
\qquad
\begin{array}c
\alga=\ddu{i-1}{i+1}i
\\[12pt]
\abacus(mmm,bbb,nbb,nbn,nnn)
\\[20pt]
\abacus(mmm,bbb,bnb,nnb,nnn)
\end{array}
\qquad
\begin{array}c
\gaal=\dw i
\\[12pt]
\abacus(mmm,bbb,bnb,nnb,nnn)
\\[20pt]
\abacus(mmm,bbb,nbb,nbn,nnn)
\end{array}
\qquad
\begin{array}c
\bega=\dl{i+1}i1
\\[12pt]
\abacus(mmm,bbb,bnb,nbn,nnn)
\\[20pt]
\abacus(mmm,bbb,nbb,nnb,nnn)
\end{array}
\]
Observe that $\albe\doms\alga,\gaal\doms\bega$ (and in fact $\alga$ and $\gaal$ are incomparable even in the usual dominance order).

\subsection{Proving \cref{main} in $B$}\label{pf3sec}

Now we consider what conditions we need on a restricted \bip $\mu$ in order to have $\dn\la\mu>1$ for some $\la$.

\begin{lemma}\label{domcond3}
Suppose $\mu$ is a restricted \bip in $B$, and $\dn\la\mu>1$ for some $\la$. Then $\nor i\mu=\nor{i+1}\mu=1$, and at least three of the decomposition numbers $\dn\albe\mu$, $\dn\alga\mu$, $\dn\gaal\mu$ and $\dn\bega\mu$ are positive. Hence $\mu\dm\dom\alga,\gaal\dom\mu$.
\end{lemma}

\begin{pf}
In the same way as in the proof of \cref{domcond}, we can apply the restriction functor $\ee i$ to show that $\mu$ has only one normal $i$-node and that if either of the decomposition numbers $\dn\albe\mu$ and $\dn\gaal\mu$ is greater than $1$ then they are both positive, while if either of $\dn\alga\mu$ and $\dn\bega\mu$ is greater than $1$ then they are both positive. But in the same way we can use $\ee{i+1}$ to show that $\nor{i+1}\mu=1$ and either $\dn\albe\mu$ and $\dn\alga\mu$ are both positive or $\dn\gaal\mu$ and $\dn\bega\mu$ are both positive. Combining these statements shows that at least three of the given decomposition numbers are positive. Now \cref{basicdecompjs}, together with the \js dominance order on the \bips $\albe,\alga,\gaal,\bega$ gives the second statement.
\end{pf}

So to prove our main result for $B$ we can focus on \bips $\mu$ in $B$ for which:
\begin{enumerate}
\item\label{isres3}
$\mu$ is restricted;
\item\label{aggadm}
$\alga,\gaal\dom\mu$;
\item\label{mdagga}
$\mu\dm\dom\alga,\gaal$; and
\item\label{allnorm3}
$\nor i\mu=\nor{i+1}\mu=1$, and $\nor j\mu=0$ for $j\neq i,i+1$.
\end{enumerate}

We begin by finding the \bips $\mu$ for which conditions (\ref{aggadm}) and (\ref{allnorm3}) are satisfied. A careful analysis of the possible abacus configurations shows that this happens in the cases given in the following table.

\begin{longtable}{|c|c|}
\hline
$\mu$&conditions
\\
\hline
$\dl{i+1}{m}1$&$l<m$
\\
$\dl{i+1}{i-1}2$&$-$
\\
$\dl{i+1}{i}2$&$-$
\\
$\dl{i+1}{l}2$&$k<l=m$
\\
$\dl{i+1}{m}2$&$l<m$
\\
$\ddu{i+1}{i+2}{m}$&$i+1<j,\ l<m$
\\
$\ddu{i+1}{k+1}{m}$&$k<l<m$
\\\hline
\end{longtable}

(Note in particular that the \bips $\mu=\dl{i+1}{i+1}2$ and $\dl{i+1}k2$ (with $i+1<k$) do not appear; these \bips satisfy (\ref{aggadm}) and (\ref{allnorm3}) if the usual dominance order is used, but $\alga\ndom\mu$ in the \js dominance order.)

We next want to check which of these \bips are restricted. To do this we use (and extend) the technique that we used for the \bip $\dl i{i-1}2$ in \cref{proving2}. Given a \bip $\mu$ in the above table, we repeatedly remove good nodes until we reach a \bip $\xi$ which has weight $0$ or $1$, or has no normal nodes. If $\xi$ has no normal nodes (but is not the empty \bip) then $\xi$ is not restricted, and so neither is $\mu$. If $\xi$ has weight $0$ then it is restricted; if $\xi$ has weight $1$, then it is restricted \iff it is not the most dominant \bip in its block.

If we determine that $\mu$ is restricted via this technique, we can also use it to construct $\mu\dm$. Suppose we have reached the restricted \bip $\xi$ of weight $0$ or $1$ by removing good nodes of residues $r_1,\dots,r_t$ in turn. The \bip $\xi\dm$ can easily be identified: if $\xi$ has weight $0$ then $\xi\dm=\xi$, while if $\xi$ has weight $1$ then (from \cite[Lemma 4.9]{mfweight}) $\xi\dm$ is the minimal \bip in the same block as $\xi$ for which $\xi\dm\doms\xi$. Having identified $\xi\dm$, we add anticogood nodes of residues $r_t,\dots,r_1$ in turn, and by \cref{dmbij} we obtain~$\mu\dm$.

\begin{eg}
Suppose $k=m$ and let $\mu=\dl{i+1}{i-1}2$, with the following \abd. (We modify the \abd by adjusting the numbers of beads in order to make runner $i-1$ the leftmost runner).
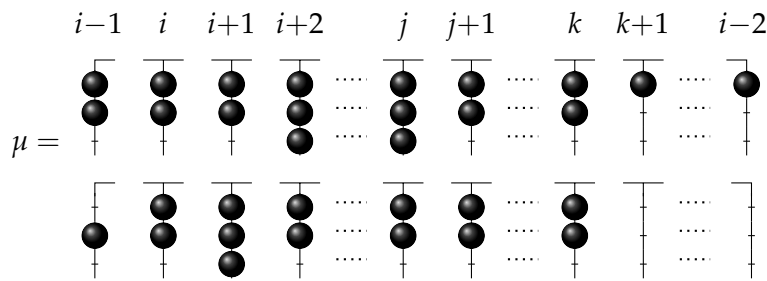
\[
\mu=
\begin{array}{cccc@{\ \ }c@{\ \ }cc@{\ \ }c@{\ \ }cc@{\ \ }c@{\ \ }c}
\mathclap{i{-}1}&\mathclap{i}&\mathclap{i{+}1}&\mathclap{i{+}2}&&\mathclap{j}&\mathclap{j{+}1}&&\mathclap{k}&\mathclap{k{+}1}&&\mathclap{i{-}2}
\\[6pt]
\abacus(l,b,b,n)&
\abacus(m,b,b,n)&
\abacus(m,b,b,n)&
\abacus(m,b,b,b)&
\abacus(h,h,h)&
\abacus(m,b,b,b)&
\abacus(m,b,b,n)&
\abacus(h,h,h)&
\abacus(m,b,b,n)&
\abacus(m,b,n,n)&
\abacus(h,h,h)&
\abacus(r,b,n,n)
\\[20pt]
\abacus(l,n,b,n)&
\abacus(m,b,b,n)&
\abacus(m,b,b,b)&
\abacus(m,b,b,n)&
\abacus(h,h,h)&
\abacus(m,b,b,n)&
\abacus(m,b,b,n)&
\abacus(h,h,h)&
\abacus(m,b,b,n)&
\abacus(m,n,n,n)&
\abacus(h,h,h)&
\abacus(r,n,n,n)
\end{array}
\]
Now we remove good nodes of residues
\[
i,\ i+1,i+1,i+2,i+2,\dots,j,j,\ j+1,\dots,k,\ i-1,i-2,\dots,k+1
\]
in turn, to obtain the \bip $\xi$ with the following \abd.
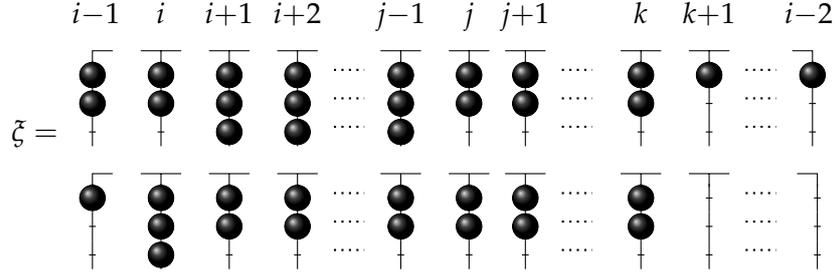
\[
\xi=
\begin{array}{cccc@{\ \ }c@{\ \ }cc@{\ \ }c@{\ \ }ccc@{\ \ }c@{\ \ }c}
\mathclap{i{-}1}&\mathclap{i}&\mathclap{i{+}1}&\mathclap{i{+}2}&&\mathclap{j{-}1}&\mathclap{j}&\mathclap{j{+}1}&&\mathclap{k}&\mathclap{k{+}1}&&\mathclap{i{-}2}
\\[6pt]
\abacus(l,b,b,n)&
\abacus(m,b,b,n)&
\abacus(m,b,b,b)&
\abacus(m,b,b,b)&
\abacus(h,h,h)&
\abacus(m,b,b,b)&
\abacus(m,b,b,n)&
\abacus(m,b,b,n)&
\abacus(h,h,h)&
\abacus(m,b,b,n)&
\abacus(m,b,n,n)&
\abacus(h,h,h)&
\abacus(r,b,n,n)
\\[20pt]
\abacus(l,b,n,n)&
\abacus(m,b,b,b)&
\abacus(m,b,b,n)&
\abacus(m,b,b,n)&
\abacus(h,h,h)&
\abacus(m,b,b,n)&
\abacus(m,b,b,n)&
\abacus(m,b,b,n)&
\abacus(h,h,h)&
\abacus(m,b,b,n)&
\abacus(m,n,n,n)&
\abacus(h,h,h)&
\abacus(r,n,n,n)
\end{array}
\]
By calculating the values $\ga_i(\xi)$ we find that $\xi$ has weight $1$. We now look at the other \bips in the same block as $\xi$, and we find that $\xi$ is restricted, and that $\xi\dm$ is obtained by moving a bead from component $2$ to component $1$ on runner $i$, and from component $1$ to component $2$ on runner $i-1$. This implies in particular that $\mu$ is restricted.
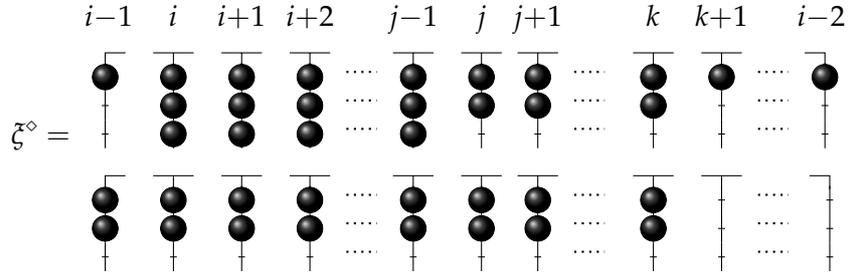
\[
\xi\dm=
\begin{array}{cccc@{\ \ }c@{\ \ }cc@{\ \ }c@{\ \ }ccc@{\ \ }c@{\ \ }c}
\mathclap{i{-}1}&\mathclap{i}&\mathclap{i{+}1}&\mathclap{i{+}2}&&\mathclap{j{-}1}&\mathclap{j}&\mathclap{j{+}1}&&\mathclap{k}&\mathclap{k{+}1}&&\mathclap{i{-}2}
\\[6pt]
\abacus(l,b,n,n)&
\abacus(m,b,b,b)&
\abacus(m,b,b,b)&
\abacus(m,b,b,b)&
\abacus(h,h,h)&
\abacus(m,b,b,b)&
\abacus(m,b,b,n)&
\abacus(m,b,b,n)&
\abacus(h,h,h)&
\abacus(m,b,b,n)&
\abacus(m,b,n,n)&
\abacus(h,h,h)&
\abacus(r,b,n,n)
\\[20pt]
\abacus(l,b,b,n)&
\abacus(m,b,b,n)&
\abacus(m,b,b,n)&
\abacus(m,b,b,n)&
\abacus(h,h,h)&
\abacus(m,b,b,n)&
\abacus(m,b,b,n)&
\abacus(m,b,b,n)&
\abacus(h,h,h)&
\abacus(m,b,b,n)&
\abacus(m,n,n,n)&
\abacus(h,h,h)&
\abacus(r,n,n,n)
\end{array}
\]
Now adding anticogood nodes of residues
\[
k+1,k+2,\dots,i-1,\ k,k-1,\dots,j+1,\ j,j,j-1,j-1,\dots,i+1,i+1,\ i
\]
in turn yields the \bip $\mu\dm=\dl{i-1}i1$.

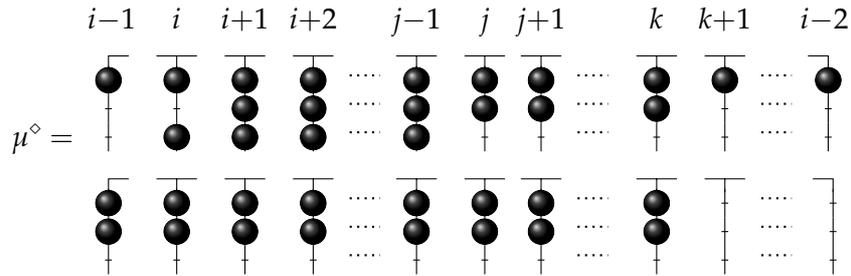
\[ 
 \mu\dm=
 \begin{array}{cccc@{\ \ }c@{\ \ }cc@{\ \ }c@{\ \ }ccc@{\ \ }c@{\ \ }c}
 \mathclap{i{-}1}&\mathclap{i}&\mathclap{i{+}1}&\mathclap{i{+}2}&&\mathclap{j{-}1}&\mathclap{j}&\mathclap{j{+}1}&&\mathclap{k}&\mathclap{k{+}1}&&\mathclap{i{-}2}
 \\[6pt]
 \abacus(l,b,n,n)&
 \abacus(m,b,n,b)&
 \abacus(m,b,b,b)&
 \abacus(m,b,b,b)&
 \abacus(h,h,h)&
 \abacus(m,b,b,b)&
 \abacus(m,b,b,n)&
 \abacus(m,b,b,n)&
 \abacus(h,h,h)&
 \abacus(m,b,b,n)&
 \abacus(m,b,n,n)&
 \abacus(h,h,h)&
 \abacus(r,b,n,n)
 \\[20pt]
 \abacus(l,b,b,n)&
\abacus(m,b,b,n)&
\abacus(m,b,b,n)&
\abacus(m,b,b,n)&
\abacus(h,h,h)&
\abacus(m,b,b,n)&
\abacus(m,b,b,n)&
\abacus(m,b,b,n)&
\abacus(h,h,h)&
\abacus(m,b,b,n)&
\abacus(m,n,n,n)&
\abacus(h,h,h)&
\abacus(r,n,n,n)
\end{array}
\]
\end{eg}

Following this procedure for all the cases in the above table yields the following list of restricted \bips $\mu$.

\needspace{3em}
\setcounter{tablecase}0
\begin{longtable}{|r|c|c|c|}
\hline
&$\mu$&conditions&$\mu\dm$
\\
\hline
\nextcase&$\dl{i+1}{i-1}2$&$k<l<m$&$\dl{i-1}{l+1}2$
\\*
\nextcase&$\dl{i+1}{i-1}2$&$k<l=m$&$\dl{i-1}{i-1}2$
\\*
\nextcase&$\dl{i+1}{i-1}2$&$k=l<m<e+i-2$&$\ddu{i-1}{i-2}{k+1}$
\\*
\rnextcase{usejsone}&$\dl{i+1}{i-1}2$&$k=l<m=e+i-2$&$\dl{i-1}{k+1}1$
\\*
\rnextcase{conj1}&$\dl{i+1}{i-1}2$&$k=l=m$&$\dl{i-1}i1$
\\
\hline
\nextcase&$\dl{i+1}i2$&$j<k<l-1$&$\ddu{l-1}l{j+1}$
\\*
\nextcase&$\dl{i+1}i2$&$j<k=l-1$&$\dl {k+1}{j+1}1$
\\*
\rnextcase{iii21}&$\dl{i+1}i2$&$j<k=l$&$\dl{i-1}{j+1}1$
\\*
\nextcase&$\dl{i+1}i2$&$j=k<l$&$\dl l{i+1}1$
\\*
\rnextcase{iii22}&$\dl{i+1}i2$&$j=k=l$&$\dl{i-1}{i+1}1$
\\
\hline
\nextcase&$\dl{i+1}l2$&$k<l-1,\ l=m$&$\ddu{l-1}li$
\\*
\nextcase&$\dl{i+1}l2$&$k=l-1,\ l=m$&$\dl{k+1}i1$
\\
\hline
\nextcase&$\dl{i+1}m2$&$k<l-1,\ l<m$&$\ddu{l-1}li$
\\*
\nextcase&$\dl{i+1}m2$&$k=l-1,\ l<m$&$\dl{k+1}i1$
\\*
\rnextcase{conj2}&$\dl{i+1}m2$&$k=l<m$&$\dl{i-1}i1$
\\
\hline
\nextcase&$\ddu{i+1}{i+2}{m}$&$i+1<j,\ k<l<m$&$\dl li2$
\\*
\nextcase&$\ddu{i+1}{i+2}{m}$&$i+1<j,\ k=l<m$&$\dl{i-1}{i-1}2$
\\
\hline
\nextcase&$\ddu{i+1}{k+1}{m}$&$k<l<m$&$\dl{i-1}{i-1}2$
\\
\hline
\end{longtable}

Our next task is to check which of these cases satisfy $\mu\dm\dom\alga,\gaal$. We find that this happens in cases \ref{conj1}, \ref{iii21}, \ref{iii22} and \ref{conj2}. (Note that in case \ref{usejsone} 
we need to bring into play the \js order: in this case $\mu\dm\dom\alga$ with the usual dominance order, but $\mu\dm\ndom\alga$ with the \js order.) Now we consider the remaining four cases.

\subsection{Cases \ref{iii21} and \ref{iii22}}

Assume $k=l$ and let $\mu=\dl{i+1}i2$. In this case we will show that $\dn\alga\mu=\dn\gaal\mu=0$, so that $\dn\la\mu\ls1$ for all $\la$ by \cref{domcond3}. To do this, we use the \js formula. First we find all \bips in the set
\[
D=\lset{\la}{\text{$\la\dom\mu$ and either $\alga\dom\la$ or $\gaal\dom\la$}}.
\]
By examining \abds, we find that
\begin{align*}
D={}&\lset{\dl{i+1}x2}{i\ls x\ls k}\cup\lset{\dw x}{j+1\ls x\ls k}
\\
&\cup\lset{\dl xx1}{i+1\ls x\ls j}\cup\{\bega,\alga,\gaal\}.
\end{align*}
The \js order on these \bips is given in \cref{case810js}.
\begin{figure}[ht]
\[
\begin{tikzpicture}[scale=.4,yscale=.8,every node/.style={fill=white,inner sep=.5pt}]
\draw(0,0)coordinate(c1);
\draw(-3,3)coordinate(c2);
\draw(-5,5)coordinate(c3a);
\draw(-7,7)coordinate(c3b);
\draw(-9,9)coordinate(c4);
\draw(-12,12)coordinate(d1);
\draw(-14,14)coordinate(d2a);
\draw(-16,16)coordinate(d2b);
\draw(-18,18)coordinate(d3);
\draw(-21,21)coordinate(c5);
\draw(-18,24)coordinate(d4);
\draw(-16,26)coordinate(d5a);
\draw(-14,28)coordinate(d5b);
\draw(-12,30)coordinate(d6);
\draw(-9,33)coordinate(c8);
\draw(-7,35)coordinate(c9a);
\draw(-5,37)coordinate(c9b);
\draw(-3,39)coordinate(c10);
\draw(0,42)coordinate(c11);
\draw(0,46)coordinate(c12);
\draw(0,21)coordinate(c6);
\draw(6,26)coordinate(c7);
\draw(c3a)--(c2)--(c1)--(c6)--(c7);
\draw(c6)--(c11)--(c12);
\draw(c11)--(c10)--(c9b);
\draw(c9a)--(c8)--(d5b);
\draw(d5a)--(d4)--(c5)--(d3)--(d2b);
\draw(d2a)--(d1)--(c4)--(c3b);
\draw[dashed](c3a)--(c3b);
\draw[dashed](c9a)--(c9b);
\draw[dashed](d2a)--(d2b);
\draw[dashed](d5a)--(d5b);
\draw(c1)node{$\dl{i+1}i2\mathrlap{=\mu}$};
\draw(c2)node{$\dl{i+1}{i+2}2$};
\draw(c4)node{$\dl{i+1}j2$};
\draw(d1)node{$\dl{i+1}{j+1}2$};
\draw(d3)node{$\dl{i+1}k2$};
\draw(c5)node{$\dl{i+1}{i+1}2$};
\draw(d4)node{$\dw k$};
\draw(d6)node{$\dw{j+1}$};
\draw(c6)node{$\dl{i+1}i1\mathrlap{=\bega}$};
\draw(c7)node{$\ddu{i-1}{i+1}i\mathrlap{=\alga}$};
\draw(c8)node{$\dl jj1$};
\draw(c10)node{$\dl{i+2}{i+2}1$};
\draw(c11)node{$\dl{i+1}{i+1}1$};
\draw(c12)node{$\dw i\mathrlap{=\gaal}$};
\end{tikzpicture}
\]
\caption{\js interval for Cases \ref{iii21} and \ref{iii22}}\label{case810js}
\end{figure}
Now we compute the decomposition numbers $\dn\la\mu$ for $\la\in D$ using the \js formula. In \cref{iii21table} we collate the needed \js coefficients $j_{\la\nu}$, and compute the bounds $J_{\la\mu}$, which give the decomposition numbers. (Note that for $\la=\dl{i+1}{i+1}1$, the bound $J_{\la\mu}$ equals $2$, but we use the fact that (by the inductive assumption, and because $\la$ is not \tricky) $\dn\la\mu\ls1$, and therefore $\dn\la\mu=1$.)

\begin{table}[ht]
\[
\begin{array}{|c|ccccc|c|c|}
\hline
\multirow{2}{*}{\raisebox{-3pt}{$\la$}}&&&\nu&&&\multirow{2}{*}{\raisebox{-3pt}{\smash{$J_{\la\mu}$}}}&\multirow{2}{*}{\raisebox{-3pt}{\smash{$\dn\la\mu$}}}
\\
&\mu&\bega&\dl{i+1}{i+2}2&\dl{i+2}{i+2}1&\dl{i+1}{i+1}1&&
\\[3pt]\hline
\bega&1&&&&&1&1
\\[3pt]
\alga&-1&1&&&&0&0
\\[3pt]
\dl{i+1}{i+2}2&1&&&&&1&1
\\[3pt]
\dl{i+1}x2\ (i+3\ls x\ls k)&(-1)^{x-i}&&(-1)^{x-i+1}&&&0&0
\\[3pt]
\dl{i+1}{i+1}2&(-1)^{k-i+1}&&(-1)^{k-i}&&&0&0
\\[3pt]
\dw x\ (k\gs x\gs j+1)&0&&0&&&0&0
\\[3pt]
\dl xx1\ (j\gs x\gs i+3)&0&&0&&&0&0
\\[3pt]
\dl{i+2}{i+2}1&0&&1&&&1&1
\\[3pt]
\dl{i+1}{i+1}1&0&1&0&1&&2&1
\\[3pt]
\gaal&-1&1&0&-1&1&0&0
\\\hline
\end{array}
\]
\caption{\js calculations for cases \ref{iii21} and \ref{iii22}}\label{iii21table}
\end{table}

We conclude that if $\la\in D$, then
\[
\dn\la\mu=
\begin{cases}
1&\text{if }\la\in\{\mu,\bega,\dl{i+1}{i+2}2,\dl{i+2}{i+2}1,\dl{i+1}{i+1}1\}
\\
0&\text{otherwise}.
\end{cases}
\]
Now the fact that $\dn\alga\mu=\dn\gaal\mu=0$, combined with \cref{domcond3}, shows that $\dn\la\mu\ls1$ for all $\la$.

\subsection{Case \ref{conj1}}\label{3case5}

Suppose $k=l=m$, and let $\mu=\dl{i+1}{i-1}2$, so that $\mu\dm=\dl{i-1}i1$. In this case we use a conjugation argument. If we let $B'$ be the block containing the \bip $(\mu\dm)'$, then $B'$ is also a block of type III, with \xr $\xi'$, and can be treated in the same way as the block $B$, with the integers $i,j,k,l,m$ replaced by $i',j',k',l',m'$, where
\begin{align*}
i'&=k_1+k_2-i-1,
\\
j'&=k_1+k_2+e-m,
\\
k'&=k_1+k_2+e-l,
\\
l'&=k_1+k_2+e-k,
\\
m'&=k_1+k_2+e-j.
\end{align*}
In the block $B'$, $(\mu\dm)'$ is the \bip $\dl{i'+1}{i'}2$, with $j'=k'=l'$. So using case \ref{iii22} dealt with above, we can deduce that $\dn\la{(\mu\dm)'}\ls1$ for all $\la$, and therefore (by \cref{conjdual}) $\dn\la\mu\ls1$ for all $\la$.

\subsection{Case \ref{conj2}}

In this case we use a conjugation argument, as we did for case \ref{conj1} above. In this case $(\mu\dm)'$ satisfies the conditions of case \ref{iii21} above, therefore we can deduce that $\dn\la\mu\ls1$ for all $\la$ in the same way.

\medskip

This deals with all cases, and so our inductive step for type III blocks is complete.

\section{Blocks of type IV}\label{type4sec}

The last situation we have to investigate is when $B$ is a block of $\hhh n$ and $i\in\mathbb{Z}$ is such that $\de_i(B)=2$ and $\de_j(B)\le0$ for every $j\nequiv i\ppmod e$. Again we will follow an inductive approach to show that \cref{main} holds for $B$.

\subsection{The abacus for a block of type IV}

Take $\xi$ the \xr of $B$ and consider its $(k_1+1,k_2-1)$-\abd following the same convention of \cref{type2sec,type3sec}. Then $\ga_{i}(\xi)=2$ and $\ga_j(\xi)\le0$ for every $j\nequiv i\ppmod e$ and we note that $\xi$ has exactly two removable $i$-nodes (again because $\xi$ has weight 0 and by \cref{36lem} it cannot have both removable and addable nodes of the same residue). We now observe that these two removable nodes occur in different components of $\xi$. Assuming the contrary we would have that $\ga_{i}(\xi)-\ga_{i-1}(\xi)=2$ but $\xi$ has weight 0 and ${\ga_{x}(\xi)-\ga_{y}(\xi)}\le1$ for every $x,y\in\zez$.

Similarly to the previous sections we draw the $(k_1+1,k_2-1)$-\abd of $\xi$ finding integers $i\ls j\ls k\ls l\ls m\ls e+i-2$ such that $\xi$ has a removable $i$-node in both components and addable nodes of residue $j+1$ and $l+1$ in component 1 and of residue $k+1$ and $m+1$ in component 2 (in particular $\kappa_1=\widebar{j+l-i+1}$, and $\kappa_2=\widebar{k+m+3-i}$). Now $Z=\{\ol i,\dots,\ol j\}\cup\{\widebar{k+1},\dots,\widebar l\}\cup\{\widebar{m+1},\dots,\widebar{i-1}\}$.

\begin{eg}
  Suppose $e=17$, $\kappa=(\ol{14},\ol5)$ and $\xi=\bp{7^4}{4^7}$, with the following abacus display.
  
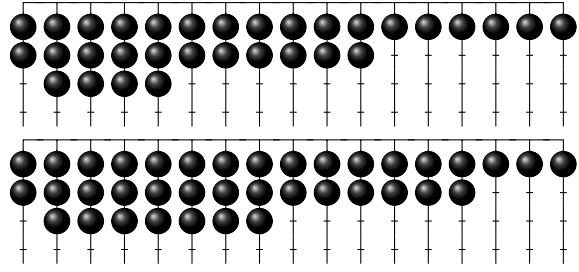
\[
  \begin{array}c
  \abacus(lmmmmmmmmmmmmmmmr,bbbbbbbbbbbbbbbbb,bbbbbbbbbbbnnnnnn,nbbbbnnnnnnnnnnnn,nnnnnnnnnnnnnnnnn)
  \\[20pt]
  \abacus(lmmmmmmmmmmmmmmmr,bbbbbbbbbbbbbbbbb,bbbbbbbbbbbbbbnnn,nbbbbbbbnnnnnnnnn,nnnnnnnnnnnnnnnnn)
  \end{array}
  \]
  We can take $(i,j,k,l,m)=(1,4,7,10,13)$ then having $Z=\{\ol0,\ol1,\ol2,\ol3,\ol4,\ol8,\ol9,\widebar{10},\widebar{14},\widebar{15},\widebar{16}\}$.
  \end{eg}
  
  \subsection{\Tricky \bips in $B$}

There are four \tricky \bips in $B$, which we label $\be_1,\be_2,\be_3,\be_4$, with abacus displays (for runners $i-1,i,i+1$) as follows.
  
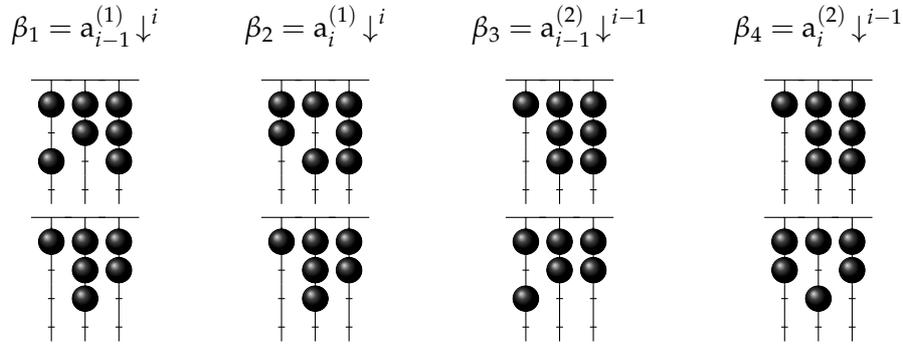
\[
\begin{array}c
\be_1=\dl i{i-1}1
\\[12pt]
\abacus(mmm,bbb,nbb,bnb,nnn)
\\[20pt]
\abacus(mmm,bbb,nbb,nbn,nnn)
\end{array}
\qquad
\begin{array}c
\be_2=\dl i{i}1
\\[12pt]
\abacus(mmm,bbb,bnb,nbb,nnn)
\\[20pt]
\abacus(mmm,bbb,nbb,nbn,nnn)
\end{array}
\qquad
\begin{array}c
\be_3=\dl{i-1}{i-1}2
\\[12pt]
\abacus(mmm,bbb,nbb,nbb,nnn)
\\[20pt]
\abacus(mmm,bbb,nbb,bnn,nnn)
\end{array}
\qquad
\begin{array}c
\be_4=\dl{i-1}i2
\\[12pt]
\abacus(mmm,bbb,nbb,nbb,nnn)
\\[20pt]
\abacus(mmm,bbb,bnb,nbn,nnn)
\end{array}
\]
In contrast to Case III, these \bips are totally ordered by the \js dominance order; in particular $\be_4\doms\be_3\doms\be_2\doms\be_1$.

\subsection{Proving \cref{main} in $B$}

Also in this situation it is useful to establish a lemma that brings out a useful property of those restricted \bips $\mu$ such that $\dn\la\mu>1$ for some $\la$. By \cref{simpledn0}(2) we note that such \bips have exactly two normal $i$-nodes. 

\begin{lemma}\label{domcond4}
  Suppose $\mu$ is a restricted \bip in $B$, and $\dn\la\mu>1$ for some $\la$. Then the decomposition numbers $\dn{\be_1}\mu,\dn{\be_2}\mu,\dn{\be_3}\mu,\dn{\be_4}\mu$ 
  are all positive; in particular $\mu\dm\dom\be_4\doms\be_3\doms\be_2\doms\be_1\dom\mu$.
\end{lemma}
  
  \begin{pf}
  Let $\la$ be a \bip such that $\dn\la\mu>1$. Then $\la=\be_x$ for some $x\in\{1,2,3,4\}$. We follow the same idea as in \cref{domcond,domcond3} with the only difference that here we need to apply the restriction functor $\ee i$ to $\spe{\be_1},\dots,\spe{\be_4}$ 'twice', namely we apply the functor $\ee i^2$ described in \cref{branchsec}. By the branching rules we find that there are \bips $\hat\be_1,\hat\be_2,\hat\be_3,\hat\be_4,\tilde\mu$ of $n-2$ such that
  \begin{align*}
    \ee i^2\spe{\be_1}&\sim(\spe{\hat\be_1})^2+(\spe{\hat\be_2})^2+(\spe{\hat\be_3})^2,
    \\
    \ee i^2\spe{\be_2}&\sim(\spe{\hat\be_1})^2+(\spe{\hat\be_2})^2+(\spe{\hat\be_4})^2,
    \\
    \ee i^2\spe{\be_3}&\sim(\spe{\hat\be_1})^2+(\spe{\hat\be_3})^2+(\spe{\hat\be_4})^2,
    \\
    \ee i^2\spe{\be_4}&\sim(\spe{\hat\be_2})^2+(\spe{\hat\be_3})^2+(\spe{\hat\be_4})^2,
    \\
    \ee i^2\smp{\mu}&\sim\smp{\tilde\mu}.
  \end{align*}
Analogously as in \cite[Proposition 2.10(4)]{mfwt3}, we find
    \[
    \dn{\be_1}\mu+\dn{\hat\be_4}{\tilde\mu}=\dn{\be_2}\mu+\dn{\hat\be_3}{\tilde\mu}=\dn{\be_3}\mu+\dn{\hat\be_2}{\tilde\mu}=\dn{\hat\be_4}\mu+\dn{\hat\be_1}{\tilde\mu}.
    \]
    and we conclude by induction that if one of the decomposition numbers $\dn{\be_x}\mu$ is greater than $1$, then they are all positive. The last sentence follows from \cref{basicdecompjs}.
  \end{pf}
  
  Summing up, it remains to consider those \bips $\mu$ in $B$ such that

\begin{enumerate}
  \item\label{isres4}
  $\mu$ is restricted;
  \item\label{bejsdm}
  $\be_1\dom\mu$;
  \item\label{mjsdbe}
  $\mu\dm\dom\be_4$; and
  \item\label{allnorm4}
  $\nor i\mu=2$, and $\nor j\mu=0$ for $j\neq i$.
  \end{enumerate}

Firstly we list all possible \bips $\mu$ which satisfy (\ref{bejsdm}) and (\ref{allnorm4}) and in which both components are $e$-restricted (which is a necessary condition for $\mu$ to be restricted, by \cref{resterest}).
\needspace{4em}
\begin{longtable}{|c|c|}
\hline
$\mu$&conditions
\\
\hline
$\dl i{m}1$&$l<m$
\\
$\dl i{i-1}2$&$-$
\\
$\dl i{m}2$&$l<m$
\\
$\ddu{i}{i+1}{m}$&$i<j,\ l<m$
\\
$\ddu{i}{k+1}{m}$&$k<l<m$
\\\hline
\end{longtable}

Now we find which \bips in this list satisfy condition (\ref{isres4}), i.e. we find the restricted \bips. To do this we use the same technique as in the previous sections. In each case where $\mu$ is restricted we find $\mu\dm$ using the procedure explained in \cref{pf3sec}.

\setcounter{tablecase}0
\needspace{3em}
\begin{longtable}{|r|c|c|c|}
\hline
&$\mu$&conditions&$\mu\dm$
\\
\hline
\nextcase&$\dl i{i-1}2$&$k<l<m$&$\dl{i-1}{l+1}2$
\\*
\nextcase&$\dl i{i-1}2$&$k<l=m$&$\dl{i-1}{i-1}2$
\\*
\nextcase&$\dl i{i-1}2$&$k=l<m<e+i-2$&$\ddu{i-1}{i-2}{k+1}$
\\*
\rnextcase{usejs1}&$\dl i{i-1}2$&$k=l<m=e+i-2$&$\dl{i-1}{k+1}1$
\\*
\rnextcase{iv12}&$\dl i{i-1}2$&$j<k=l=m$&$\dl{i-1}{j+1}1$
\\*
\rnextcase{iv13}&$\dl i{i-1}2$&$j=k=l=m$&$\dl{i-1}i1$
\\
\hline
\nextcase&$\dl im2$&$j<k<l-1,\ l<m$&$\ddu{l-1}l{j+1}$
\\*
\nextcase&$\dl im2$&$j<k=l-1,\ l<m$&$\dl{k+1}{j+1}{1}$
\\*
\rnextcase{iv31}&$\dl im2$&$j<k=l<m$&$\dl{i-1}{j+1}1$
\\*
\nextcase&$\dl im2$&$j=k<l<m$&$\dl{l}{i}1$
\\*
\rnextcase{iv324}&$\dl im2$&$j=k=l<m$&$\dl{i-1}{i}1$
\\
\hline
\nextcase&$\ddu{i}{i+1}{m}$&$i<j,\ k<l<m$&$\dl li2$
\\*
\nextcase&$\ddu{i}{i+1}{m}$&$i<j,\ k=l<m$&$\dl{i-1}{i-1}2$
\\
\hline
\nextcase&$\ddu{i}{k+1}{m}$&$k<l<m$&$\dl{i-1}{i-1}{2}$
\\
\hline
\end{longtable}

Finally we cut our list with condition (\ref{mjsdbe}). It turns out that there are four cases left: \ref{iv12}, \ref{iv13}, \ref{iv31}, and \ref{iv324}. (Note that case \ref{usejs1} is the only one where we need to use the \js dominance order in condition (\ref{mjsdbe}) instead of the usual dominance order).

\subsection{Cases \ref{iv12} and \ref{iv324}}\label{4case5/13}

To address cases \ref{iv12} and \ref{iv324}, we use a conjugation argument similar to that used for cases \ref{conj1} and \ref{conj2} in \cref{type3sec}. If $\mu$ satisfies the conditions of case \ref{iv12}, then $\mu\dm=\dl{i-1}{j+1}1$. So $(\mu\dm)'$ satisfies the conditions of case \ref{iv324} for the block with \xr $\xi'$, and conversely. So, by \cref{conjdual}, showing that $\dn\la\mu\ls1$ for all $\la$ in case \ref{iv12} is equivalent to doing the same in case \ref{iv324}.

In fact, we will need to consider both cases simultaneously: we will obtain a partial result in case \ref{iv324}, and use this in our analysis of case \ref{iv12} to obtain a complete result. So we assume first that $\mu$ satisfies the conditions of case \ref{iv324}. Our objective here will be to prove the following.

\begin{lemma}\label{case13claim}
Suppose $\mu$ satisfies the conditions of case \ref{iv324}. Then $\dn{\be_1}\mu=1$, while $\dn{\be_2}\mu\gs\dn{\be_3}\mu\gs1$.
\end{lemma}

\begin{figure}[hb]
\[
\begin{tikzpicture}[scale=.4,yscale=.8,xscale=-1,every node/.style={fill=white,inner sep=1pt}]
\draw(0,-5)coordinate(c1);
\draw(-5,1)coordinate(c2);
\draw(-5,3)coordinate(c3a);
\draw(-5,6)coordinate(c3b);
\draw(-5,8)coordinate(c4);
\draw(-5,12)coordinate(d1);
\draw(-5,16)coordinate(c5);
\draw(-5,18)coordinate(d5a);
\draw(-5,21)coordinate(d5b);
\draw(-5,23)coordinate(c10);
\draw(5,29)coordinate(c12);
\draw(5,25)coordinate(c6);
\draw(5,-1)coordinate(c7);
\draw(5,3)coordinate(cc);
\draw(5,5)coordinate(c7a);
\draw(5,8)coordinate(c7b);
\draw(5,10)coordinate(c16);
\draw(5,14)coordinate(c17);
\draw(5,16)coordinate(c17a);
\draw(5,19)coordinate(c17b);
\draw(5,21)coordinate(c18);
\draw(5,34)coordinate(c13);
\draw(16,14)coordinate(c19);
\draw(c3a)--(c2)--(c1);
\draw(c1)--(c7)--(cc)--(c7a);
\draw(c7b)--(c16)--(c17)--(c17a);
\draw(c17b)--(c18)--(c6)--(c12)--(c13);
\draw(c7)--(c19)--(c13);
\draw(d5b)--(c10)--(c12);
\draw(d5a)--(c5)--(d1)--(c4)--(c3b);
\draw[dashed](c3a)--(c3b);
\draw[dashed](c7a)--(c7b);
\draw[dashed](c17a)--(c17b);
\draw[dashed](d5a)--(d5b);

\draw(c1)node{$\dl i{m}2\mathrlap{=\mu}$};
\draw(c2)node{$\dl i{i+1}2$};
\draw(c4)node{$\dl ij2$};
\draw(d1)node{$\dl i{i}2$};
\draw(c5)node{$\dl jj1$};
\draw(c10)node{$\dl{i+1}{i+1}1$};
\draw(c7)node{$\ddu{i}{i-1}{m}$};
\draw(cc)node{$\ddu{i}{e+i-2}{m}$};
\draw(c16)node{$\ddu{i}{m+1}{m}$};
\draw(c17)node{$\dl i{m}{1}$};
\draw(c18)node{$\dl i{e+i-2}{1}$};
\draw(c6)node{$\dl i{i-1}1\mathrlap{=\be_1}$};
\draw(c12)node{$\dl ii1\mathrlap{=\be_2}$};
\draw(c13)node{$\dl {i-1}{i-1}2\mathrlap{=\be_3}$};
\draw(c19)node{$\dl{i-1}{m}{2}$};
\end{tikzpicture}
\]
\caption{\js interval for case \ref{iv324}}\label{case13js}
\end{figure}

To prove \cref{case13claim}, we compute decomposition numbers using the \js formula. For this, we need to consider the set of \bips
\[
D=\lset{\la}{\text{$\be_3\dom\la\dom\mu$}}.
\]
By considering possible \abds, we find that
\begin{align*}
D={}&\lset{\dl ix2}{i\ls x\ls j}\cup\lset{\dl xx1}{i+1\ls x\ls j}\cup\lset{\ddu{i}{x}{m}}{m+1\ls x\ls e+i-2}\\
&\cup\lset{\dl i{x}{1}}{m\ls x\ls e+i-2}\cup\left\{\mu,\dl {i-1}m2,\ddu i{i-1}m,\be_1,\be_2,\be_3\right\},
\end{align*}
and that the \js order on these \bips is given in \cref{case13js}.
Now, as in the cases in \cref{type3sec} we collect in a table the integers $j_{\la\nu}$ and $J_{\la{\mu}}$ needed to calculate decomposition numbers $\dn\la{\mu}$. We obtain the results in \cref{iv324table}, from which \cref{case13claim} follows.



\begin{table}[!hb]
{\small
\[
\begin{array}{|c|cccccc|c|c|}
\hline
\multirow{2}{*}{\raisebox{-3pt}{$\la$}}&&&\nu&&&&\multirow{2}{*}{\raisebox{-3pt}{\smash{$J_{\la{\mu}}$}}}&\multirow{2}{*}{\raisebox{-3pt}{\smash{$\dn\la{\mu}$}}}
\\
&\mu&\dl i{i+1}2&\dl{i+1}{i+1}1&\ddu{i}{i-1}{m}&\be_1&\be_2&&
\\[3pt]\hline
\dl i{i+1}2&1&&&&&&1&1
\\[3pt]
\dl ix2\ (i+2\ls x\ls j)&(-1)^{x-i+1}&(-1)^{x-i}&&&&&0&0
\\[3pt]
\dl ii2&(-1)^{j-i}&(-1)^{j-i+1}&&&&&0&0
\\[3pt]
\dl xx1\ (j\gs x\gs i+2)&0&0&&&&&0&0
\\[3pt]
\dl{i+1}{i+1}1&0&1&&&&&1&1
\\[3pt]
\ddu{i}{i-1}{m}&1&&&&&&1&1
\\[3pt]
\dl{i-1}{m}{2}&-1&&1&&&&0&0
\\[3pt]
\begin{matrix}\ddu{i}{x}{m}\\(e+i-2\gs x\gs m+1)\end{matrix}&(-1)^{e+x-i+1}&&(-1)^{e+x-i}&&&&0&0
\\[3pt]
\dl i{m}{1}&(-1)^{e+m-i+1}&&(-1)^{e+m-i}&&&&0&0
\\[3pt]
\begin{matrix}\dl i{x}{1}\\(m+1\ls x\ls e+i-2)\end{matrix}&0&&0&&&&0&0
\\[3pt]
\be_1&0&&&1&&&1&1
\\[3pt]
\be_2&0&0&1&0&1&&2&
\\[3pt]
\be_3&0&0&-1&0&1&1&\dn{\be_2}\mu&
\\\hline
\end{array}
\]
}
\caption{\js calculations for case \ref{iv324}}\label{iv324table}
\end{table}
Having proved \cref{case13claim}, we now analyse case \ref{iv12} fully. Assume $j<k=l=m$ and $\mu=\dl i{i-1}{2}$. First we use \cref{conjdual} to derive some consequences of \cref{case13claim}: since conjugation reverses the (\js) dominance order on \bips, the conjugates of the \bips $\be_1,\be_2,\be_3,\be_4$ from case \ref{iv324} are the \bips $\be_4,\be_3,\be_2,\be_1$ in case \ref{iv12}. So \cref{conjdual,case13claim} imply that in case \ref{iv12}
\[
\dn{\be_4}\mu=1\qquad\text{and}\qquad1\ls\dn{\be_2}\mu\ls\dn{\be_3}\mu.
\]
We define $\tau$ be the \bip $\dl{i-1}{i+1}2$, and we want to consider the decomposition numbers $\dn\la\mu$ for \bips $\la$ in the set
\[
D=\lset{\la}{\text{$\tau\dom\la\dom\mu$}}.
\]
We can compute
\[
D=\lset{\dl ix2}{i-1\ls x\ls k}\cup\lset{\dw x}{j+1\ls x\ls k}\cup\lset{\dl xx1}{i+1\ls x\ls j}\cup\left\{\be_1,\be_2,\be_3,\be_4,\tau\right\}
\]
and the \js order on the \bips in $D$ is given in \cref{case5js}.
\begin{figure}[!b]
\[
\begin{tikzpicture}[scale=.4,yscale=.8,xscale=.9,every node/.style={fill=white,inner sep=1pt}]
\draw(0,-1)coordinate(c1);
\draw(-5,2)coordinate(c2);
\draw(-7.5,3.5)coordinate(c3a);
\draw(-10,5)coordinate(c3b);
\draw(-12.5,6.5)coordinate(c4);
\draw(-17.5,9.5)coordinate(d1);
\draw(-20,11)coordinate(d2a);
\draw(-22.5,12.5)coordinate(d2b);
\draw(-25,14)coordinate(d3);
\draw(-30,17)coordinate(c5);
\draw(-25,20)coordinate(d4);
\draw(-22.5,21.5)coordinate(d5a);
\draw(-20,23)coordinate(d5b);
\draw(-17.5,24.5)coordinate(d6);
\draw(-12.5,27.5)coordinate(c8);
\draw(-10,29)coordinate(c9a);
\draw(-7.5,30.5)coordinate(c9b);
\draw(-5,32)coordinate(c10);
\draw(0,35)coordinate(c12);
\draw(0,16)coordinate(c6);
\draw(0,39)coordinate(c13);
\draw(0,43)coordinate(c14);
\draw(0,47)coordinate(c15);
\draw(c3a)--(c2)--(c1)--(c6);
\draw(c12)--(c10)--(c9b);
\draw(c9a)--(c8)--(d5b);
\draw(d5a)--(d4)--(c5)--(d3)--(d2b);
\draw(d2a)--(d1)--(c4)--(c3b);
\draw(c6)--(c12)--(c13)--(c14)--(c15);
\draw[dashed](c3a)--(c3b);
\draw[dashed](c9a)--(c9b);
\draw[dashed](d2a)--(d2b);
\draw[dashed](d5a)--(d5b);

\draw(c1)node{$\dl i{i-1}2\mathrlap{=\mu}$};
\draw(c2)node{$\dl i{i+1}2$};
\draw(c4)node{$\dl ij2$};
\draw(d1)node{$\dl i{j+1}2$};
\draw(d3)node{$\dl ik2$};
\draw(c5)node{$\dl i{i}2$};
\draw(d4)node{$\dw k$};
\draw(d6)node{$\dw{j+1}$};
\draw(c6)node{$\dl i{i-1}1\mathrlap{=\be_1}$};
\draw(c8)node{$\dl jj1$};
\draw(c10)node{$\dl{i+1}{i+1}1$};
\draw(c12)node{$\dl ii1\mathrlap{=\be_2}$};
\draw(c13)node{$\dl {i-1}{i-1}2\mathrlap{=\be_3}$};
\draw(c14)node{$\dl {i-1}i2\mathrlap{=\be_4}$};
\draw(c15)node{$\dl{i-1}{i+1}2\mathrlap{=\tau}$};
\end{tikzpicture}
\]
\caption{\js interval for case \ref{iv12}}\label{case5js}
\end{figure}
Now we apply the \js formula to compute decomposition numbers $\dn\la\mu$ for $\be_2\doms\la\dom\mu$. These are given in \cref{case5table}.
\begin{table}[t]
\[
\begin{array}{|c|cc|c|c|}
\hline
\multirow{2}{*}{\raisebox{-3pt}{$\la$}}&\nu&&\multirow{2}{*}{\raisebox{-3pt}{\smash{$J_{\la{\mu}}$}}}&\multirow{2}{*}{\raisebox{-3pt}{\smash{$\dn\la{\mu}$}}}
\\
&\mu&\dl i{i+1}2&&
\\[3pt]\hline
\dl i{i+1}2&1&&1&1
\\[3pt]
\dl ix2\ (i+2\ls x\ls k)&(-1)^{x-i+1}&(-1)^{x-i}&0&0
\\[3pt]
\dl ii2&(-1)^{k-i}&(-1)^{k-i+1}&0&0
\\[3pt]
\dw x\ (k\gs x\gs j+1)&0&0&0&0
\\[3pt]
\dl xx1\ (j\gs x\gs i+2)&0&0&0&0
\\[3pt]
\dl{i+1}{i+1}1&0&1&1&1
\\[3pt]
\be_1&1&&1&1
\\\hline
\end{array}
\]
\caption{\js calculations for case \ref{iv12}}\label{case5table}
\end{table}

We also use the \js formula to estimate $\dn\tau\mu$: we obtain
\[
J_{\tau\mu}=-\dn{\dl i{i+1}2}\mu+\dn{\dl{i+1}{i+1}1}\mu-\dn{\be_3}\mu+\dn{\be_4}\mu=1-\dn{\be_3}\mu,
\]
forcing $\dn{\be_3}\mu=1$. Therefore we have $\dn{\be_x}\mu=1$ for $x=1,2,3,4$.

\subsection{Case \ref{iv31}}

The analysis of case \ref{iv31} is very similar to the cases in \cref{4case5/13}. We begin by proving an exact analogue of \cref{case13claim} for case \ref{iv31}. The calculation is very similar to the calculation used to prove \cref{case13claim}: if $\be_3\dom\la\doms\mu$, then we find that
\[
J_{\la\mu}=\begin{cases}
  1 & \text{if}\ \la\in\{\dl i{i+1}{2},\ \dl{i+1}{i+1}{1},\ \ddu{i}{i-1}{m},\ \be_1\},
  \\
  2 & \text{if}\ \la=\be_2,
  \\
  \dn{\be_2}\mu & \text{if}\ \la=\be_3,
  \\
  0 &  \text{otherwise}.
\end{cases}
\]
(Note that if $i=j$ then a slight modification is necessary: the \bip $\dw{i+1}$ replaces $\dl{i+1}{i+1}{1}$ in the above formula.)

Having established the analogue of \cref{case13claim}, we apply duality: case \ref{iv31} is self-dual, so we deduce that
\[
\dn{\be_1}\mu=1,\qquad\dn{\be_2}\mu=\dn{\be_3}\mu\gs1,\qquad\dn{\be_4}\mu=1
\]
whenever $\mu$ satisfies the conditions of case \ref{iv31}. We complete case \ref{iv31} in the same way as case \ref{iv12} above, letting $\tau$ be the \bip $\dl{i-1}{i+1}2$. As in case \ref{iv12}, we find that $J_{\tau\mu}=1-\dn{\be_3}\mu$, forcing $\dn{\be_3}\mu=\nolinebreak1$.

\subsection{Case \ref{iv13}}

This is the last case we have to deal with. Assume that $j=k=l=m$ and $\mu=\dl i{i-1}{2}$. In this case we also begin by proving a direct analogue of \cref{case13claim}. This is proved as in the previous cases: the \bips $\la$ for which $\be_3\dom\la\dom\mu$ are shown in \cref{case6js} (ignoring the red portion of the diagram), and calculating \js coefficients as in previous cases we find that if $\be_3\dom\la\doms\mu$, then
\[
J_{\la\mu}=\begin{cases}
  1 & \text{if}\ \la\in\{\dl i{i+1}{2},\ \dl{i+1}{i+1}{1}, \be_1\},
  \\
  2 & \text{if}\ \la=\be_2,
  \\
  \dn{\be_2}\mu & \text{if}\ \la=\be_3,
  \\
  0 &  \text{otherwise}.
\end{cases}
\]
So an analogue of \cref{case13claim} follows. As with case \ref{iv31}, this case is self-dual, so we deduce that
\[
\dn{\be_1}\mu=1,\qquad\dn{\be_2}\mu=\dn{\be_3}\mu\gs1,\qquad\dn{\be_4}\mu=1
\]
whenever $\mu$ satisfies the conditions of case \ref{iv13}.

Now assume that $j<e+i-2$. With this assumption we complete case \ref{iv13} in the same way as in the preceding cases. Let $\tau$ be the \bip $\dl{e+i-2}i2$. \cref{case6js} (now including the red part) shows the \js order on \bips $\la$ for which $\tau\dom\la\dom\mu$.

\begin{figure}[!t]
\[
\begin{tikzpicture}[scale=.4,yscale=.8,every node/.style={fill=white,inner sep=1pt}]
\draw(0,0)coordinate(c1);
\draw(-4,4)coordinate(c2);
\draw(-6,6)coordinate(c2a);
\draw(-8,8)coordinate(c2b);
\draw(-10,10)coordinate(c3);
\draw(-14,14)coordinate(c4);
\draw(-10,18)coordinate(c5);
\draw(-8,20)coordinate(c5a);
\draw(-6,22)coordinate(c5b);
\draw(-4,24)coordinate(c6);
\draw(0,14)coordinate(c8);
\draw(0,28)coordinate(c9);
\draw(17,51)coordinate(c2a');
\draw(17,49)coordinate(c2b');
\draw(10,36)coordinate(c6');
\draw(10,41)coordinate(c7');
\draw(0,38)coordinate(c8');
\draw(0,33)coordinate(c9');
\draw(10,31)coordinate(c10');
\draw(c1)--(c2)--(c2a);
\draw(c2b)--(c3)--(c4)--(c5)--(c5a);
\draw(c5b)--(c6)--(c9);
\draw(c1)--(c8)--(c9)--(c9');
\draw[red](c9')--(c8');
\draw[red](c6')--(c9');
\draw[red](c7')--(c8');
\draw[red](c7')--(c6')--(c10')--(c1);
\draw[dashed](c2a)--(c2b);
\draw[dashed](c5a)--(c5b);
\draw(c1)node{$\dl i{i-1}2\mathrlap{=\mu}$};
\draw(c2)node{$\dl i{i+1}{2}$};
\draw(c3)node{$\dl i{j}{2}$};
\draw(c4)node{$\dl i{i}{2}$};
\draw(c5)node{$\dl{j}{j}{1}$};
\draw(c6)node{$\dl{i+1}{i+1}{1}$};
\draw(c8)node{$\dl i{i-1}{1}\mathrlap{=\be_1}$};
\draw(c9)node{$\dl i{i}{1}\mathrlap{=\be_2}$};
\draw(c6')node[text=red]{$\dl{e+i-2}{e+i-2}{2}$};
\draw(c7')node[text=red]{${\dl{e+i-2}{i}{2}\mathrlap{=\tau}}$};
\draw[red](c8')node{$\mathllap{\be_4=}\dl{i-1}{i}{2}$};
\draw(c9')node{$\mathllap{\be_3=}\dl{i-1}{i-1}{2}$};
\draw(c10')node[text=red]{$\dl{e+i-2}{i-1}{2}$};
\end{tikzpicture}
\]
\caption{\js interval for case \ref{iv13}}\label{case6js}
\end{figure}

As in the preceding cases, we find that $J_{\tau\mu}=1-\dn{\be_3}\mu$, forcing $\dn{\be_3}\mu=\dn{\be_2}\mu=1$, so that we are done when $j<e+i-2$. By duality, case \ref{iv13} is also complete when $i<j$.

%
%
%
%

Then the only situation left to be checked is when $i=j=k=l=m=e+i-2$. In this case $e=2$ and $B$ is the block of $\hhh 5$ which contains eight Specht modules labelled by the following \bips.

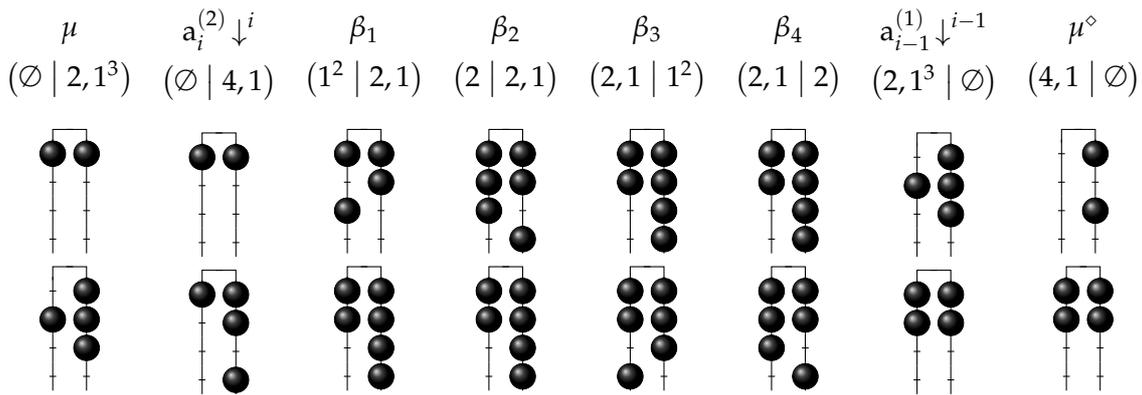
\[
\begin{array}c
\mu
\\[5pt]
\bp{\emptyset}{2,1^3}
\\[12pt]
\abacus(lr,bb,nn,nn,nn)
\\[20pt]
\abacus(lr,nb,bb,nb,nn)
\end{array}
\begin{array}c
\dl i{i}{2}
\\[5pt]
\bp{\emptyset}{4,1}
\\[12pt]
\abacus(lr,bb,nn,nn,nn)
\\[20pt]
\abacus(lr,bb,nb,nn,nb)
\end{array}
\begin{array}c
\be_1
\\[5pt]
\bp{1^2}{2,1}
\\[12pt]
\abacus(lr,bb,nb,bn,nn)
\\[20pt]
\abacus(lr,bb,bb,nb,nb)
\end{array}
\begin{array}c
\be_2
\\[5pt]
\bp{2}{2,1}
\\[12pt]
\abacus(lr,bb,bb,bn,nb)
\\[20pt]
\abacus(lr,bb,bb,nb,nb)
\end{array}
\begin{array}c
\be_3
\\[5pt]
\bp{2,1}{1^2}
\\[12pt]
\abacus(lr,bb,bb,nb,nb)
\\[20pt]
\abacus(lr,bb,bb,nb,bn)
\end{array}  
\begin{array}c
\be_4
\\[5pt]
\bp{2,1}{2}
\\[12pt]
\abacus(lr,bb,bb,nb,nb)
\\[20pt]
\abacus(lr,bb,bb,bn,nb)
\end{array}
\begin{array}c
  \dl{i-1}{i-1}{1}
  \\[5pt]
  \bp{2,1^3}{\emptyset}
  \\[12pt]
  \abacus(lr,nb,bb,nb,nn)
  \\[20pt]
  \abacus(lr,bb,bb,nn,nn)
  \end{array}
\begin{array}c
  \mu\dm
  \\[5pt]
  \bp{4,1}{\emptyset}
  \\[12pt]
  \abacus(lr,nb,nn,nb,nn)
  \\[20pt]
  \abacus(lr,bb,bb,nn,nn)
\end{array}
\]and two simple modules labelled by $\mu$ and $\be_1$. It is easy to solve this case by hand. We use a dimension argument here. First we note that the two Specht modules $\spe{\be_1}$ and $\spe{\be_2}$ have the same dimension. From our calculations above we know that
\[
\dn{\be_1}\mu=1,
\]
while \cref{basicdecompjs} and a very simple calculation with the \js formula show that
\[
\dn{\be_1}{\be_1}=\dn{\be_2}{\be_1}=1,
\]
forcing $\dn{\be_2}\mu=1$.

This completes the proof of Case IV, and of our main theorem.

\end{document}